%&amstex            
\input amstex\documentstyle{amsppt}  
\pagewidth{12.5cm}\pageheight{19cm}\magnification\magstep1
\topmatter
\title Hecke algebras and involutions in Coxeter groups\endtitle
\author G. Lusztig and D. A. Vogan, Jr.\endauthor
\address{Department of Mathematics, M.I.T., Cambridge, MA 02139}\endaddress
\thanks{Supported in part by National Science Foundation grants DMS-1303060 and DMS-0967272.}\endthanks   
\endtopmatter   
\document

\define\coker{\text{\rm coker}}

\define\mpb{\medpagebreak}

\define\bF{\bar F}

\define\hR{\hat R}

\define\si{\sim}

\define\qua{\quad}

\define\lb{\linebreak}

\define\op{\oplus}
   
\define\part{\partial}

\define\imp{\implies}

\define\n{\notin}
\define\iy{\infty}
\define\m{\mapsto}
\define\do{\dots}

\define\lra{\leftrightarrow}

\define\sub{\subset}

\define\ti{\tilde}
\define\nl{\newline}
\redefine\i{^{-1}}

\define\un{\underline}
\define\ov{\overline}
\define\ot{\otimes}

\define\Hom{\text{\rm Hom}}

\define\tr{\text{\rm tr}}

\define\uuc{\un{\un{\car}}}
\define\a{\alpha}
\redefine\b{\beta}
\redefine\c{\chi}

\redefine\d{\delta}
\define\e{\epsilon}

\redefine\o{\omega}

\define\ph{\phi}

\define\s{\sigma}
\redefine\t{\tau}

\define\x{\xi}

\define\Th{\Theta}

\define\Ph{\Phi}
\define\Ps{\Psi}

\redefine\aa{\bold a}

\define\HH{\bold H}
\define\II{\bold I}

\define\NN{\bold N}

\define\RR{\bold R}

\define\ZZ{\bold Z}

\define\ca{\Cal A}
\define\cb{\Cal B}

\define\ch{\Cal H}

\define\ck{\Cal K}

\define\cm{\Cal M}

\define\car{\Cal R}

\define\cv{\Cal V}

\define\cz{\Cal Z}

\define\cy{\Cal Y}

\define\ff{\frak f}

\define\fh{\frak h}

\define\fC{\frak C}

\define\tf{\ti f}

\define\tw{\ti w}
\define\tz{\ti z}

\define\ty{\ti y}
\define\tA{\ti A}

\define\tC{\ti C}

\define\tM{\ti M}

\define\tP{\ti P}

\define\sh{\sharp}

\define\bul{\bullet}

\define\che{\check}

\define\cir{\circ}

\define\EW{EW}
\define\KL{KL}
\define\QG{L1}
\define\HEC{L2}
\define\INV{L3}
\define\LV{LV}
\define\SO{S}

\head Introduction\endhead
\subhead 0.1\endsubhead
Let $W$ be a finitely generated Coxeter group with a fixed involutive automorphism $w\m w^*$ which leaves 
stable the set of simple reflections. An element $w\in W$ is said to be a $*$-twisted involution if 
$w\i=w^*$. Let $\II=\{w\in W;w\i=w^*\}$ be the set of $*$-twisted involutions of $W$. Let $\ca'=\ZZ[v,v\i]$ 
where $v$ is an indeterminate. In \cite{\LV} we have defined (geometrically) an action of the Hecke algebra 
of $W$ (with parameter $v^2$) on the free $\ca'$-module $\cm$ with basis $\{a_w;w\in\II\}$, assuming that 
$W$ is a Weyl group. In \cite{\INV} a definition of the Hecke algebra action on $\cm$ was given in a purely 
algebraic way, without assumption on $W$. The purpose of this paper is to give a more conceptual approach to
the definition of the Hecke algebra action on $\cm$, based on the theory of Soergel bimodules \cite{\SO} and
on the recent results of Elias and Williamson \cite{\EW} in that theory. 

In this paper we interpret $\cm$ as a (modified) Grothendieck group associated to the category of Soergel
bimodules corresponding to $W$ and to a $2$-periodic functor of this category to itself, defined using $*$
and by switching left and right multiplication in a bimodule. The action of the Hecke algebra appears quite 
naturally in this interpretation; however, we must find a way to compute explicitly the action of a 
generator $T_s+1$ of the Hecke algebra ($s$ is a simple reflection) 
on a basis element $a_w$ of $\cm$ so that we recover the formulas of \cite{\LV}, \cite{\INV}. 
The formula has four cases depending on whether $sw$ is equal to $ws^*$ or not and on whether the length
of $sw$ is smaller or larger than that of $w$. In each case, $(T_s+1)a_w$ is a linear combination 
$c'a_{w'}+c''a_{w''}$ of two basis elements $a_{w'},a_{w''}$ where one of $w',w''$ is equal to $w$, the other
is $sw$ or $sws^*$ and the length of $w'$ is smaller than that of $w''$. We cannot prove the formulas 
directly. Instead we compute directly the coefficient $c'$ and then observe that if $c'$ is known, then 
$c''$ is automatically known from the fact that we have a Hecke algebra action. The computation of $c'$ 
occupies Sections 4 and 5 (see Theorem 5.2). It has two cases (depending on whether $sw'$ is equal to 
$w's^*$ or not). The two cases require quite different proofs.

As an application of Theorem 6.2 (which is essentially a corollary of Theorem 5.2) we outline a proof (6.3) 
of a positivity conjecture (9.12 in \cite{\INV}) stating that, if $y,w\in\II$ and $\d\in\{1,-1\}$, then the 
polynomial $P^\s_{y,w}$ introduced in \cite{\INV} (and earlier in \cite{\LV} in the case of Weyl groups) 
satisfies $(P_{y,w}(u)+\d P^\s_{y,w}(u))/2\in\NN[u]$ where $P_{y,w}$ is the polynomial introduced in 
\cite{\KL}. This is a refinement of the statement \cite{\EW} that $P_{y,w}(u)\in\NN[u]$ which holds for any 
$y,w\in W$.
In \S7 we show as another application of our results that $\cm$ admits a filtration by Hecke algebra
submodules whose subquotients are indexed by the two-sided cells of $W$. Under a boundedness assumption we
show that the Hecke algebra acts on such a subquotient by something resembling a $W$-graph.

\head Contents\endhead

1. $2$-periodic functors.

2. A review of Soergel modules.

3. The $\HH$-module $\cm$.

4. Some exact sequences.

5. Trace computations.

6. Applications.

7. The $\HH$-module $\cm_c$.

\head 1. $2$-periodic functors\endhead
\subhead 1.1\endsubhead
In this section we review some results from \cite{\QG, \S11}.

Let $k$ be a field of characteristic zero. Let $\fC$ be a $k$-linear category, that is a category in which 
the space of morphisms between any two objects has a given $k$-vector space structure such that composition 
of morphisms is bilinear and such that finite direct sums exist. A functor from one $k$-linear category to 
another is said to be $k$-linear if it respects the $k$-vector space structures.

Let $\ck(\fC)$ be the Grothendieck group of $\fC$ that is, the free abelian group generated by symbols $[A]$
for each $A\in\fC$ (up to isomorphism) with relations $[A\op B]=[A|+[B]$ for any $A,B\in\fC$.
A $k$-linear functor $M\m M^\sh$, $\fC@>>>\fC$ is said to be {\it $2$-periodic} if $M\m(M^\sh)^\sh$ is the 
identity functor $\fC@>>>\fC$. Assuming that such a functor is given we define a new $k$-linear category 
$\fC_\sh$ as follows. The objects of $\fC_\sh$ are pairs $(A,\ph)$ where $A\in\fC$ and $\ph:A^\sh@>>>A$ is an
isomorphism in $\fC$ such that the composition $(A^\sh)^\sh@>\ph^\sh>>A^\sh@>\ph>>A$ is the identity map of 
$A$. Let $(A,\ph)$, $(A',\ph')$ be two objects of $\fC_\sh$. We define a $k$-linear map 
$\Hom_{\fC}(A,A')@>>>\Hom_{\fC}(A,A')$ by $f\m f^!:=\ph'f^\sh\ph\i$. Note that $(f^!)^!=f$. By definition,
$\Hom_{\fC_\sh}((A,\ph),(A',\ph'))=\{f\in\Hom_{\fC}(A,A');f=f^!\}$, a $k$-vector space. 
The direct sum of two objects $(A,\ph)$, $(A',\ph')$ is $(A\op A',\ph\op\ph')$. Clearly, if 
$(A,\ph)\in\fC_\sh$,
then $(A,-\ph)\in\fC_\sh$. An object $(A,\ph)$ of $\fC_\sh$ is said to be {\it traceless} if there exists an 
object $B$ of $\fC$ and an isomorphism $A\cong B\op B^\sh$ under which $\ph$ corresponds to an isomorphism 
$B^\sh\op B@>\si>>B\op B^\sh$ which carries the first (resp. second) summand of $B^\sh\op B$ onto the second
(resp. first) summand of $B\op B^\sh$.

Let $\ck_\sh(\fC)$ be the quotient of $\ck(\fC_\sh)$ by the subgroup $\ck^0(\fC_\sh)$ generated by the 
elements $[B,\ph]$ where $(B,\ph)$ is any traceless object of $\fC_\sh$. We show that:
$$[A,-\ph]=-[A,\ph]\text{ for any }(A,\ph)\in\fC_\sh.\tag a$$
Indeed, if we define $\ph':A^\sh\op A@>>>A\op A^\sh$ by $(x,y)\m(y,x)$ and
$\t:A\op A@>>>A\op A^\sh$ by $(x,y)\m(x+y,\ph\i(x)-\ph\i(y))$, then $\t$ defines an isomorphism of
$(A,\ph)\op(A,-\ph)$ with the traceless object $(A\op A^\sh,\ph')$.

\head 2. A review of Soergel modules\endhead
\subhead 2.1\endsubhead
In this section we review some results of Soergel \cite{\SO} and of Elias-Williamson \cite{\EW}.

Recall that $W$ is a Coxeter group. The canonical set of generators (assumed to be finite) is denoted by 
$S$. Let $x\m l(x)$ be the length function on $W$ and let $\le$ be the Bruhat order on $W$. Let 
$\fh$ be a reflection representation of $W$ over the real numbers $\RR$, as in \cite{\EW}; for any $s\in S$ 
we fix a linear form $\a_s:\fh@>>>\RR$ whose kernel is equal to the fixed point set of $s:\fh@>>>\fh$. Let 
$R$ be the algebra of polynomial functions $\fh@>>>\RR$ with the $\ZZ$-grading in which linear functions
$\fh@>>>\RR$ have degree $2$. Note that $W$ acts naturally on $R$; we write this action as $w:r\m{}^wr$ and 
for $s\in S$ we set $R^s=\{r\in R;{}^sr=r\}$, a subalgebra of $R$. Let $R^{>0}=\{r\in R;r(0)=0\}$.
Let $\hR$ be the completion of $R$ with respect to the maximal ideal $R^{>0}$. 

Let $\car$ be the category whose objects are $\ZZ$-graded $(R,R)$-bimodules in which for $M,M'\in\car$,
$\Hom_\car(M,M')$ is the space of  homomorphisms of $(R,R)$-bimodules $M@>>>M'$ compatible with the
$\ZZ$-gradings. For $M\in\car$ and $n\in\ZZ$, the shift $M[n]$ is the object of $\car$ equal in degree $i$ to
$M$ in degree $i+n$. For $M,M'$ in $\car$ we set $MM'=M\ot_RM'$; this is naturally an object of $\car$. For 
$M,M'$ in $\car$ we set 
$$M'{}^M=\op_{n\in\ZZ}\Hom_\car(M,M'[n]),$$
viewed as an object of $\car$ with $(rf)(m)=f(rm)$, $(fr)(m)=f(mr)$ for $m\in M,f\in M'{}^M,r\in R$.
For any $M\in\car$ we set $\un{M}=M/MR^{>0}=M\ot_R\RR$ where $\RR$ is identified with $R/R^{>0}$. We view 
$\un{M}$ as a $\ZZ$-graded $\RR$-vector space. For any $M\in\car$ we set $\hat M=M\ot_R\hR$, viewed as a
$\ZZ$-graded right $\hR$-module.

For $s\in S$ let $B_s=R\ot_{R^s}R[1]\in\car$.
More generally, for any $x\in W$, Soergel \cite{\SO, 6.16} shows that there is an object $B_x$ of $\car$
(unique up to isomorphism) such that $B_x$ is an indecomposable direct summand of $B_{s_1}B_{s_2}\do B_{s_q}$
for some/any reduced expression $w=s_1s_2\do s_q$ ($s_i\in S$) and such that $B_x$ is not a direct summand 
of $B_{s'_1}B_{s'_2}\do B_{s'_p}$ whenever $s'_1,\do,s'_p\in S,p<q$.
Let $\tC$ be the full subcategory of $\car$ whose objects are isomorphic to finite direct sums of shifts of 
objects of the form $B_x$ for various $x\in W$. Let $C$ be the full subcategory of $\car$ whose objects are
isomorphic to finite direct sums of objects of the form $B_x$ for various $x\in W$.
From \cite{\SO} it follows that for $M,M'\in\tC$ we have $MM'\in\tC$.

(In the case where $W$ is a Weyl group of a reductive group $G$, $C$ can be thought of as the category of
semisimple $G$-equivariant perverse sheaves on the product $\cb^2$ of two copies of the flag manifold and 
$\tC$ can be thought of as the category whose objects are complexes of sheaves on $\cb^2$ which are 
(non-canonically) direct sums of objects of semisimple $G$-equivariant perverse sheaves with shifts. Then
$M,M'\m MM'$ corresponds to convolution of complexes of sheaves.)

For any $x\in W$ let $R_x$ be the object of $\car$ such that $R_x=R$ as a left $R$-module and such that for
$m\in R_x,r\in R$ we have $mr=({}^xr)m$. The following result appears in \cite{\SO, 6.15}:

(a) {\it For any $M\in\tC$, $R_x^M$ is a finitely generated graded free right $R$-module; hence 
$\dim_\RR\un{R_x^M}<\iy$.}
\nl
Note that $\un{R_x^{M[n]}}_i=\un{R_x^M}_{i-n}$ for any $i,n\in\ZZ$.

(In the case where $W$ is a Weyl group of a reductive group $G$ then 
$\un{R_x^M}_i$ can be thought of as the dual of a stalk of a cohomology sheaf of a complex of sheaves
on $\cb^2$ at a point in the $G$-orbit on $\cb^2$ corresponding to $x$.)

Let $t\in\Hom_\car(B_s[-1],R_s)=(R_s^{B_s})_1$ be the unique element such that $t(1\ot\a_s+1\ot\a_s)=0$, 
$t(1\ot1)=1$. The image of $t$ in $\un{R_s^{B_s}}_1$ is an $\RR$-basis of this one-dimensional $\RR$-vector 
space. Hence we have canonically $\un{R_s^{B_s}}_1=\RR$.

\subhead 2.2\endsubhead
Let $x\in W$. From \cite{\EW} it follows that $\Hom_\car(B_x,B_x)=\RR$ 
and from \cite{\SO, 6.16} it follows that $\dim\un{R_x^{B_x}}_{l(x)}=1$. Thus 
$\un{R_x^{B_x}}_{l(x)}\ot_\RR B_x$ is an object of $C$ isomorphic to $B_x$ and
defined up to unique isomorphism
(even though $B_x$ was defined only up to non-unique isomorphism). From now on we will use the notation 
$B_x$ for this new object. 

It satisfies 
$$\un{R_x^{B_x}}_{l(x)}=\RR.$$
When $x=s\in S$, this agrees with the earlier description of $B_s$.

\subhead 2.3\endsubhead
Let $x,x'\in W$, $x\ne x'$. From \cite{\EW} it follows that $\Hom_\car(B_{x'},B_x)=0$. This, together with 
the equality $\Hom_\car(B_x,B_x)=\RR$ implies that the objects $B_x$ are simple in $C$. Conversely, it is 
clear that any simple object of $C$ is isomorphic to some $B_x$.

\subhead 2.4\endsubhead
Let $\ca=\ZZ[u,u\i]$ where $u$ is an indeterminate.
Let $\HH$ be the free $\ca$-module with basis $T_w,w\in W$. It is known that there is a unique associative
$\ca$-algebra structure on $\HH$ such that $T_wT_{w'}=T_{ww'}$ whenever $l(ww')=l(w)+l(w')$ and
$T_s^2=u^2T_1+(u^2-1)T_s$ for $s\in S$. Note that $T_1$ is a unit element. Let $\{c_w;w\in W\}$ be the 
$\ca$-basis of $\HH$ which in \cite{\KL} was denoted by $\{C'_w;w\in W\}$. Recall that 
$$c_w=u^{-l(w)}\sum_{y\le w}P_{y,w}(u^2)T_y\tag a$$
where $P_{y,w}=1$ if $y=w$ and $P_{y,w}$ is a polynomial of degree $\le (l(w)-l(y)-1)/2$ if $y<w$.
We regard $\ck(\tC)$ as an $\ca$-module by $u^n[M]=[M[-n]]$ for $M\in\tC,n\in\ZZ$.
Note that $\ck(\tC)$ is an associative $\ca$-algebra with product defined by $[M][M']=[MM']$ for 
$M\in\tC,M'\in\tC$. From \cite{\SO, 1.10, 5.3} we see that 

(b) {\it the assignment $M\m\sum_{y\in W,i\in\ZZ}\dim\un{R_y^M)}_iu^{-i}T_y$ defines an $\ca$-algebra 
isomorphism $\c:\ck(\tC)@>\si>>\HH$.}
\nl
From \cite{\EW, Theorem 1.1} it follows that
$$\c(B_x)=c_x.\tag c$$

\head 3. The $\HH$-module $\cm$\endhead
\subhead 3.1\endsubhead 
In this section we preserve the setup of Section 2.
Recall that $w\m w^*$ is an involutive automorphism $W@>\si>>W$ leaving $S$ stable.
We can assume that there exists an involutive $\RR$-linear map $\fh@>>>\fh$ (denoted again by $x\m x^*$) 
which satisfies $(wx)^*=w^*x^*$ for $w\in W,x\in\fh$ and satisfies $\a_{s^*}=(\a_s)^*$ for $s\in S$. We fix 
such a linear map. It induces a ring involution $R@>>>R$ denoted again by $r\m r^*$.
For $M\in\car$ let $M^\sh$ be the object of $\car$ which is equal to $M$ as a graded $\RR$-vector space, but
left (resp. right) multiplication by $r\in R$ on $M^\sh$ equals right (resp. left) multiplication by $r^*$ 
on $M$. Clearly, $(M^\sh)^\sh=M$. If $f:M_1@>>>M_2$ is a morphism in $\car$ then $f$ can be also viewed as a 
morphism $M^\sh_1@>>>M^\sh_2$ in $\car$.
Note that $M\m M^\sh$ is an $\RR$-linear, $2$-periodic functor $\car@>>>\car$. Hence $\car_\sh$ is well 
defined, see 1.1. 

If $M_1,M_2\in\car$ then we have an obvious identification $M^\sh_1M^\sh_2=(M_2M_1)^\sh$ as objects in 
$\car$ (it is given by $x_1\ot x_2\m x_2\ot x_1$).

Let $s\in S$. The $\RR$-linear isomorphism $\o_s:B_{s^*}[-1]@>\si>>B_s[-1]$ given by
$x\ot_{R^{s^*}}y\m y^*\ot_{R^s}x^*$ for $x,y\in R$ 
can be viewed as an isomorphism $B^\sh_{s^*}[-1]@>\si>>B_s[-1]$ in $\car$ or as an isomorphism 
$B^\sh_{s^*}@>\si>>B_s$ in $\car$.

Now let $x\in W$ and let $s_1s_2\do s_k$ be a reduced expression for $x$. Since $B_x$ is an indecomposable 
direct summand of $B_{s_1}B_{s_2}\do B_{s_k}$ (and $k$ is minimal with this property) we see that $B^\sh_x$ 
is an indecomposable direct summand of 
$$(B_{s_1}B_{s_2}\do B_{s_k})^\sh=B^\sh_{s_k}\do B^\sh_{s_2}B^\sh_{s_1}\cong 
B_{s^*_k}\do B_{s^*_2}B_{s^*_1}$$
(and $k$ is minimal with this property) hence by \cite{\SO, 6.16} we have 
$$B^\sh_x\cong B_{(x^*)\i}.\tag a$$. 
(We use that $s^*_k\do s^*_2s^*_1$ is a reduced expression for $(x^*)\i$.)
In particular we have $B^\sh_x\in C$. It follows that $M\in C\imp M^\sh\in C$ and $M\in\tC\imp M^\sh\in\tC$.
Note that $M\m M^\sh$ are $\RR$-linear, $2$-periodic functors $C@>>>C$ and $\tC@>>>\tC$. 
Hence $C_\sh,\tC_\sh$ are defined as in 1.1 and $\ck_\sh(C)$, $\ck_\sh(\tC)$ 
are well defined abelian groups. 

\subhead 3.2\endsubhead 
Recall that $\II=\{y\in W;y\i=y^*\}$. Let $x\in W$. 
We define $\ff_x:R_x^\sh@>>>R_{(x^*)\i}$ by $r\m\ff_x(r)=({}^{x\i}r)^*$. This is an isomorphism in $\car$.

Now assume that $x\in\II$; then $\ff_x:R_x^\sh@>>>R_x$ is given by $r\m\ff_x(r)={}^x(r^*)$
and $(R_x,\ff_x)\in\car_\sh$; thus $(R_x[i],\ff_x[i])\in\car_\sh$ for any $i\in\ZZ$.
Hence, if $(M,\ph)\in\tC_\sh$ and $i\in\ZZ$, then $f\m f^!$, 
$\Hom_\car(M,R_x[i])@>>>\Hom_\car(M,R_x[i])$ is defined as in 1.1.
Taking direct sum over $i\in\ZZ$ we obtain a map $f\m f^!$, $R_x^M@>>>R_x^M$ such that $(f^!)^!=f$. 
(We always write $R_x^M$ instead of $(R_x)^M$.)
From the definitions, for $f\in R_x^M,r\in R$ we have $(fr)^!=r^*f^!$, $(rf)^!=f^!r^*$.
Since for $r\in R,b\in R_x$ we have $rb=b{}^{x\i}r$ we see that $R^{>0}R_x=R_xR^{>0}$ so that
$R^{>0}(R_x^M)=(R_x^M)R^{>0}$; we see that $f\m f^!$ induces an $\RR$-linear (involutive) map
$\un{R_x^M}@>>>\un{R_x^M}$ and (for any $i$) an $\RR$-linear involutive map $\un{R_x^M}_i@>>>\un{R_x^M}_i$ 
denoted by $\cy_{x,\ph,i}^M$. Let 
$$\e^x_i(M,\ph)=\tr_\RR(\cy_{x,\ph,i}^M,\un{R_x^M})\in\ZZ.$$
We now take $M=B_x$ (still assuming $x\in\II$ so that $(B_x,\ph)\in\tC_\sh$ for some $\ph$). 
Then $\un{R_x^{B_x}}_{l(x)}=\RR$ hence $\e^x_{l(x)}(B_x,\ph)=\pm1$. We can normalize $\ph:B_x^\sh@>>>B_x$ 
uniquely so that $\e^x_{l(x)}(B_x,\ph)=1$. We shall denote this normalized $\ph$ by $\ph_x$.

Due to 2.3, we can apply \cite{\QG, 11.1.8} to $C,\sh$; we see that 

(a) $\ck_\sh(C)$ is a free abelian group with basis $\{[B_x,\ph_x];x\in\II\}$.

\subhead 3.3 \endsubhead 
Let $\ca'=\ZZ[v,v\i]$ where $v$ is an indeterminate. We view $\ca=\ZZ[u,u\i]$ as a subring of $\ca'$ by 
setting $u=v^2$. Note that $\ck_\sh(\tC)$
 can be viewed as an $\ca'$-module with $v^n[M,\ph]=[M[-n],\ph]$ for 
$(M,\ph)\in\tC_\sh$, $n\in\ZZ$. We show:

(a) {\it The map $q:\ca'\ot\ck_\sh(C)@>>>\ck_\sh(\tC)$, $v^n\ot[M,\ph]\m[M[-n],\ph]$ is an isomorphism of
$\ca'$-modules.}
\nl
The map $q$ is clearly well defined. To prove that it is surjective
we shall use the functors $M\m\t_{\le i}M$ from $\tC$ to $\tC$ (resp. $M\m\ch^iM$
from $\tC$ to $C$) defined in \cite{\EW, 6.2}. (Here $i\in\ZZ$.) These define in an obvious way functors 
$\tC_\sh@>>>\tC_\sh$ (resp. $\tC_\sh@>>>C_\sh$) denoted again by $\t_{\le i}$ (resp. $\ch^i$). Let 
$(M,\ph)\in\tC_\ph$. From the definition we have an exact sequence in $\tC$ (with morphisms in $\tC_\sh$)
$$0@>>>\t_{\le i-1}M@>e>>\t_{\le i}M@>e'>>\ch^iM[-i]@>>>0$$
which is split but the splitting is not necessarily given by morphisms in $\tC_\sh$. Thus there exist 
morphisms 
$$\t_{\le i-1}M@<f<<\t_{\le i}M@<f'<<\ch^iM[-i]$$
in $\tC$ such that $e'f'=1$, $fe=1$, $f'e'+ef=1$.
Now $f^!,f'{}^!$ are defined as in 1.1 and, since $e^!=e, e'{}^!=e!$ (notation of 1.1), we have
$e'f'{}^!=1$, $f^!e=1$, $f'{}^!e'+ef^!=1$ hence setting $\tf=(f+f^!)/2$, $\tf'=(f'+f'{}^!)/2$, we have
$e'\tf'=1$, $\tf e=1$, $\tf'e'+e\tf=1$ and $\tf^!=\tf$, $\tf'{}^!=\tf'$. Thus we obtain a new splitting of 
the exact sequence above which is given by morphisms in $\tC_\sh$. It follows that
$$(\t_{\le i}M,\ph)\cong(\t_{\le i-1}M,\ph)\op(\ch^iM[-i],\ph)$$
 in $\tC_\sh$ (the maps $\ph$ are induced
by $M^\sh@>>>M$). Hence $[\t_{\le i}M,\ph]=[\t_{\le i-1}M,\ph]+[\ch^iM[-i],\ph]$ in $\ck_\sh(\tC)$.
Since $[M,\ph]=[\t_{\le i}M,\ph]$ for $i\gg0$ and $0=[\t_{\le i}M,\ph]$ for $-i\gg0$ we deduce that
$[M,\ph]=\sum_i[\ch^iM[-i],\ph]$. This proves the surjectivity of $q$.

We define $\ck(\tC_\sh)@>>>\ca'\ot\ck(C_\sh)$ by $[M,\ph]\m\sum_{n\in\ZZ}v^{-n}[\ch^nM,\ph_n]$ where 
$\ph_n$ is induced by $\ph$. This clearly induces a homomorphism $q':\ck_\sh(\tC)@>>>\ca'\ot\ck_\sh(C)$ which
satisfies $q'q=1$. It follows that $q$ is injective, completing the proof of (a).

\subhead 3.4\endsubhead
Using 3.2(a), 3.3(a), we see that:

(a) $\ck_\sh(\tC)$ is a free $\ca'$-module with basis $\{[B_x,\ph_x];x\in\II\}$, (notation of 3.2).

\subhead 3.5\endsubhead
Let $\cm$ be the free $\ca'$-module with basis $\{a_x;x\in\II\}$. For any $(M,\ph)\in\tC_\sh$ and any 
$y\in\II$ we set
$$\e^y(M,\ph)=\sum_{i\in\ZZ}\e^y_i(M,\ph)v^{-i}\in\ca'.$$ 
The homomorphism $\ck(\tC_\sh)@>>>\cm$, 
$$[M,\ph]\m\sum_{y\in\II}\e^y(M,\ph)a_y$$ 
clearly factors through an $\ca'$-module homomorphism 
$$\c':\ck_\sh(\tC)@>>>\cm.\tag a$$
We show:
$$\c'\text{ is an isomorphism}.\tag b$$
For $x\in\II$ let $\tA_x=\c'([B_x,\ph_x])$.
We can write $\tA_x=\sum_{y\in\II}f_{y,x}a_y$ where $f_{y,x}\in\ca'$ are zero for all but finitely many $y$.
In view of 3.4(a), to prove (b) it is enough to show:

(c) {\it Let $y\in\II$. If $y\not\le x$ then $f_{y,x}=0$. If $y\le x$ then $f_{y,x}=v^{-l(x)}\tP_{y,x}(u)$ 
where $\tP_{y,x}=1$ if $y=x$ and $\tP_{y,x}$ is a polynomial with integer coefficients of degree 
$\le(l(x)-l(y)-1)/2$ if $y<x$.}
\nl
Assume that $f_{y,x}\ne0$. Then for some $i$ we have $\e^y_i(B_x,\ph_x)\ne0$ hence $\un{R_y^{B_x}}\ne0$.
Using 2.4(b),(c) we deduce that the coefficient of $T_y$ in $c_x$ is nonzero; thus we have $y\le x$, as
required. Next we assume that $y\le x$. We have
$v^{l(x)}f_{y,x}=\sum_i\e^y_i(B_x,\ph_x)v^{-i+l(x)}$ 
hence it is enough to show that 

$\e^y_i(B_x,\ph_x)\ne0$ implies $-i+l(x)\in2\NN$ and $-i+l(x)\le l(x)-l(y)$ with strict inequality unless 
$x=y$. 
\nl
Now $\e^y_i(B_x,\ph_x)\ne0$ implies $\un{R_y^{B_x}}_i\ne0$. Hence it is enough to show that 

$\un{R_y^{B_x}}_i\ne0$ implies $-i+l(x)\in2\NN$ and $-i+l(x)\le l(x)-l(y)$ with strict inequality unless 
$x=y$. 
\nl
By 2.4(a),(b),(c) we have 
$$\sum_{j\in\ZZ}\dim\un{R_y^{B_x}}_ju^{-j+l(x)}=P_{y,x}(u^2)$$
and the desired result follows from the properties of $P_{y,x}$ (see 2.4(a)). This proves (c) hence also (b).

Next we note that for $y\in\II,y\le x$ and $\d\in\{1,-1\}$ the following holds:
$$(P_{y,x}(u)+\d\tP_{y,x}(u))/2\in\NN[u].\tag d$$
We have 
$$P_{y,x}(u)+\d\tP_{y,x}(u)=\sum_{j\in\ZZ}\dim\un{R_y^{B_x}}_jv^{-j+l(x)}
+\d\sum_{j\in\ZZ}\e^y_j(B_x,\ph_x)v^{-j+l(x)}$$
hence it is enough to show that 
$$\dim\un{R_y^{B_x}}_j+\d\e^y_j(B_x,\ph_x)\in2\NN$$
This follows from the fact that for an involutive automorphism $\t$ of a real vector space $V$ we
have $\dim(V)+\d\tr(\t,V)\in2\NN$.

\subhead 3.6\endsubhead
For any $M\in\tC$ we define a functor $F_M:\tC_\sh@>>>\tC_\sh$ by 
$$(M',\ph)\m(MM'M^\sh,\ph')$$ 
where $\ph':(MM'M^\sh)^\sh=MM'{}^\sh M^\sh@>>>MM'M^\sh$ is given by 
$$m_1\ot m'\ot m_2\m m_2\ot\ph(m')\ot m_1.$$
Note that $F_M$ induces an $\ca'$-linear map $\ck(\tC_\sh)@>>>\ck(\tC_\sh)$ which clearly maps
$\ck^0(\tC_\sh)$ into itself hence it induces an $\ca'$-linear map $\bF_M:\ck_\sh(\tC)@>>>\ck_\sh(\tC)$. If 
$M_1,M_2\in\tC$ we have $F_{M_1M_2}=F_{M_1}F_{M_2}$ hence $\bF_{M_1M_2}=\bF_{M_1}\bF_{M_2}$; moreover for 
any $(M,\ph)\in\tC_\sh$ we have 
$$F_{M_1\op M_2}(M,\ph)=((M_1\op M_2)M(M_1^\sh\op M_2^\sh),\ph')
=F_{M_1}(M,\ph)\op F_{M_2}(M,\ph)\op(\tM,\ti\ph)$$
(for a suitable $\ph'$)
where $\tM=M_2MM_1^\sh\op M_1MM_2^\sh$ and $\ti\ph:\tM^\sh@>>>\tM$ are such that $(\tM,\ti\ph)$ is a 
traceless object of $\tC_\sh$. It follows that $\bF_{M_1\op M_2}=\bF_{M_1}+\bF_{M_2}$. We see that 
$[M]\m\bF_M$ makes $\ck_\sh(\tC)$ into a (left) $\ck(\tC)$-module. 
From the definitions, for any $M\in\tC,(M',\ph)\in\tC_\sh,n\in\ZZ$ we have 
$F_{M[n]}(M',\ph)=F_M(M'[2n],\ph)$.
Hence for $h\in\ck(\tC)$, $h'\in \ck_\sh(\tC)$, $n\in\ZZ$ we have $(u^nh)h'=v^{2n}(hh')=u^n(hh')$. Via the 
isomorphism $\c:\ck(\tC)@>\si>>\HH$ in 2.4(b) and the isomorphism $\c':\ck_\sh(\tC)@>\si>>\cm$ in 
3.5(a),(b), $\cm$ becomes a (left) $\HH$-module (with $u\in\HH$ acting on $\cm$ as multiplication by 
$u=v^2$).

\head 4. Some exact sequences\endhead
\subhead 4.1\endsubhead
In this section we fix $s\in S$ and we write $\a$ instead of $\a_s$ so that $\a^*=\a_{s^*}$. Let 
$R^{s^*,>0}=R^{s^*}\cap R^{>0}$. Let $\un{\un{R}}=R/R^{s^*,>0}R$, a $\ZZ$-graded $\RR$-algebra which can be 
naturally identified with $\RR[\a^*]/(\a^{*2})$ (it is zero except in degree $0$ and $2$). Let $\uuc$ be the
category whose objects are $\ZZ$-graded right $\un{\un{R}}$-modules. For any $M'\in\car$ we write 
$\un{\un{M'}}=M'/M'R^{s^*,>0}=M'\ot_{R^{s^*}}\RR$ where $\RR=R^{s^*}/R^{s^*,>0}$ is viewed as a
$R^{s^*}$-algebra in the obvious way. Note that $\un{\un{M'}}$ is naturally an object of $\uuc$.

\subhead 4.2\endsubhead
For any $M\in\car$ we write $R.M$ (resp. $M.R$) instead of $R\ot_{R^s}M\in\car$ (resp. 
$M\ot_{R^{s^*}}R\in\car$); 
for $r\in R,m\in M$ we write $r.m$ (resp. $m.r$) instead of $r\ot m\in R.M$ (resp. 
$m\ot r\in M.R$). Note that any element of $R.M$ (resp. $M.R$) can be written uniquely in the form 
$\sum_{i\in\{0,1\}}\a^i.m_i$ (resp. $\sum_{i\in\{0,1\}}m_i\a^{*i}$) where $m_i\in M$. 

For $M,N\in\car$ let ${}'\hom(M,N)$ (resp. $\hom'(M,N)$) be the set of maps $M@>>>N$ which are homomorphisms
of $(R^s,R)$-bimodules (resp. $(R,R^{s^*})$-bimodules) and are compatible with the $\ZZ$-gradings; let 
$${}'\hom^\bul(M,N)=\op_{i\in\ZZ}{}'\hom(M,N[i]),\qua \hom'{}^\bul(M,N)=\op_{i\in\ZZ}\hom'(M,N[i]).$$
The statements (i)-(ii) below are easily verified.

(i) There is a unique group isomorphism 

${}'\hom^\bul(M,N)@>\si>>N^{R.M}$ (resp. $\hom'{}^\bul(M,N)@>\si>>N^{M.R}$), 
\nl
$f\m F$, such that for $m\in M$ we have 

$F(1.m)=f(m)$, $F(\a.m)=\a f(m)$ (resp. $F(m.1)=f(m)$, $F(m.\a^*)=f(m)\a^*$);
\nl
this is in fact an isomorphism in $\car$,
provided that ${}'\hom^\bul(M,N)$ (resp. $\hom'{}^\bul(M,N)$) is viewed as an object of $\car$ with 
$(rf)(m)=r(f(m))$, $(fr)(m)=(f(m))r$ for $m\in M,r\in R$ and $f\in{}'\hom^\bul(M,N)$ (resp. 
$f\in\hom'{}^\bul(M,N)$).

(ii) The map 

$f\m G$, $G(m)=\a.f(m)+1.f(\a m)$ (resp. \lb $G(m)=f(m).\a^*+f(m\a^*).1$)
\nl
is an isomorphism 

${}''\hom^\bul(M,N[-2])@>\si>>(R.N)^M$ (resp. $\hom''{}^\bul(M,N[-2])@>>>(N.R)^M$
\nl
in $\car$, provided 
that ${}'\hom^\bul(M,N)$ (resp. $\hom'{}^\bul(M,N)$) is viewed as an object of $\car$ with $(rf)(m)=f(rm)$, 
$(fr)(m)=f(mr)$ for $m\in M,r\in R$ and $f\in{}'\hom^\bul(M,N)$ (resp. $f\in\hom'{}^\bul(M,N)$).

Combining (i),(ii) we see that:

(iii) The map $F\m G$, 

$G(m)=\a.F(1.m)+1.F(1.\a m)$ (resp. $G(m)=F(m.1).\a^*+F(m\a^*.1).1$)
\nl
is an isomorphism 

$(N[-2])^{R.M}@>\si>>(R.N)^M$ (resp. $(N[-2])^{M.R}@>\si>>(N.R)^M$)
\nl
of $(R^s,R)$-bimodules (resp. of $(R,R^{s^*})$-bimodules).
\nl
(We use that the two $(R,R)$-bimodule structures on ${}'\hom^\bul(M,N)$ described in (i),(ii) restrict to the
same $(R^s,R)$-bimodule structure and that the two $(R,R)$-bimodule structures on $\hom'{}^\bul(M,N)$ 
described in (i),(ii) restrict to the same $(R,R^{s^*})$-bimodule structure.)

\subhead 4.3\endsubhead
For any $M'\in\car$ we write $R.M'.R$ instead of $R\ot_{R^s}M'\ot_{R^{s^*}}R\in\car$; for $r,r'$ in $R$ and
$m'\in M'$ we write $r.m'.r'$ instead of $r\ot m'\ot r'\in R.M'.R$. Note that any element $\x\in R.M'.R$ can
be written uniquely in the form $\sum_{i,j\in\{0,1\}}\a^i.\x_{ij}.a^{*j}$ where $\x_{ij}\in M'$.

For $M,N\in\car$ let $\hom(M,N)$ be the set of maps $M@>>>N$ which are homomorphisms of
$(R^s,R^{s^*})$-bimodules and which are compatible with the $\ZZ$-gradings; let 
$\hom^\bul(M,N)=\op_{i\in\ZZ}\hom(M,N[i])$. 

We view $\hom^\bul(M,N)$ as an object of $\car$ in two ways: for $f\in\hom^\bul(M,N),r\in R,m\in M$ we set 
either

(a) $(rf)(m)=r(f(m))$, $(fr)(m)=(f(m))r$;

(b) or $(rf)(m)=f(rm)$, $(fr)(m)=f(mr)$.
\nl
The statements (i),(ii) below are easily verified.

(i) There is a unique group isomorphism $\hom^\bul(M,N)@>\si>>N^{R.M.R}$ in $\car$, $f\m F$ such that for
any $m\in M$ we have $F(1.m.1)=f(m)$, $F(\a.m.1)=\a f(m)$, $F(1.m.\a^*)=f(m)\a^*$,
$F(\a.m.\a^*)=\a f(m)\a^*$; this is in fact an isomorphism in $\car$ provided that $\hom^\bul(M,N)$ is viewed
as an object of $\car$ as in (a).

(ii) The map $f\m G$, 

$G(m)=1.f(\a m\a^*).1+\a.f(m\a^*).1+1.f(\a m).\a^*+\a.f(m).\a^*$
\nl
is an isomorphism $\hom^\bul(M,N[-4])@>\si>>(R.N.R)^M$ in $\car$ provided that $\hom^\bul(M,N)$ is viewed
as an object of $\car$ as in (b).

Combining (i),(ii) we see that:

(iii) The map $F\m G$, 

$G(m)=\a.F(1.m.1).\a^*+\a.F(1.m\a^*.1).1+1.F(1.\a m.1).\a^*+1.F(1.\a m\a^*.1).1$
\nl
is an isomorphism $(N[-4])^{R.M.R}@>\si>>(R.N.R)^M$ of $(R^s,R^{s^*})$-bimodules.
\nl
(We use that the two $(R,R)$-bimodule structures on $\hom^\bul(M,N)$ described in (a),(b) restrict to the
same $(R^s,R^{s^*})$-bimodule structure.)

\subhead 4.4\endsubhead
Let $M\in\tC$ and let $\o\in W$. We define an exact sequence
$$0@>>>R_\o^M@>c>>R_\o^{R.M}@>d>>R_{s\o}^M[2]\tag a$$
as follows. We identify $R_\o^{R.M}={}'\hom^\bul(M,R_\o)$ as objects
of $\car$ as in 4.2(i); then $c$ is the obvious inclusion $R_\o^M\sub{}'\hom^\bul(M,R_\o)$ and
$d:{}'\hom^\bul(M,R_\o)@>>>R_{s\o}^M[2]$ is given by $f\m f'$, where $f'(m)={}^s(f(\a m)-\a f(m))$. (The
kernel of $d$ is clearly $R_\o^M$.) Now (a) induces sequences
$$0@>>>\un{R_\o^M}@>>>\un{R_\o^{R.M}}@>>>\un{R_{s\o}^M[2]}@>>>0,\tag b$$
$$0@>>>\un{\un{R_\o^M}}@>>>\un{\un{R_\o^{R.M}}}@>>>\un{\un{R_{s\o}^M[2]}}@>>>0.\tag c$$
We state the following result.

(d) {\it If $l(\o)<l(s\o)$, then the sequences (b),(c) are exact.}
\nl
For (b) this is implicit in the proof in \cite{\SO, Proposition 5.7, Corollary 5.16} of the fact that, under
the assumption of (d), the alternating sum of dimensions of the terms of (b) is zero (in each degree).
The statement for (c) can be reduced to that for (b) as follows.
The $\un{\un{R}}$-modules in (c) are free of finite rank (we use 2.1(a)) and the kernel and cokernel of 
right multiplication by $\a^*$ in these $\un{\un{R}}$-modules form sequences which can be identified with 
the sequence (b) which are already known to be exact; it follows that the sequence (c) is exact.

Next we define an exact sequence
$$0@>>>R_\o^M@>c'>>R_\o^{M.R}@>d'>>R_{\o s^*}^M[2]\tag e$$
as follows. We identify $R_\o^{M.R}=\hom'{}^\bul(M,R_\o)$ as objects
of $\car$ as in 4.2(i); then $c'$ is the obvious inclusion $R_\o^M\sub\hom'{}^\bul(M,R_\o)$ and
$d':\hom'{}^\bul(M,R_\o)@>>>R_{\o s^*}^M[2]$ is given by $f\m f'$, where 
$f'(m)=f(m\a^*)-f(m)\a^*$ (the product $f(m)\a^*$ is computed in the right $R$-module structure of $R_\o$).
(The kernel of $d'$ is clearly $R_\o^M$.) Now (e) induces sequences
$$0@>>>\un{R_\o^M}@>>>\un{R_\o^{M.R}}@>>>\un{R_{\o s^*}^M[2]}@>>>0,\tag f$$
$$0@>>>\un{\un{R_\o^M}}@>>>\un{\un{R_\o^{M.R}}}@>>>\un{\un{R_{\o s^*}^M[2]}}@>>>0.\tag g$$

We now state the following result.

(h) {\it If $l(\o)<l(\o s^*)$, then the sequences (f),(g) are exact.}
\nl
For any $x\in W$ we have an isomorphism  
$$R_x^M@>\si>>R_{x\i}^{M^\sh}\tag i$$ 
as $\RR$-vector spaces (not in $\car$) given by $f\m\tf$ where $(\tf)(m)={}^{x\i}(f(m))$ for any $m\in M$
(we identify $M,M^\sh$ as sets). It carries $R_x^MR^{>0}$ onto $R_{x\i}^{M^\sh}R^{>0}$ hence it induces an 
isomorphism $\un{R_x^M}@>\si>>\un{R_{x\i}^{M^\sh}}$ of graded $\RR$-vector spaces. Applying an isomorphism 
like (i) to each term of the sequence (f) we get a sequence
$$0@>>>\un{R_{\o\i}^{M^\sh}}@>>>\un{R_{\o\i}^{R\ot_{R_{s^*}}M^\sh}}@>>>\un{R_{s^*\o\i}^{M^\sh}[2]}@>>>0;$$
(we use that $(M.R)^\sh=R\ot_{R_{s^*}}M^\sh$). This sequence is a special case of the sequence (b) (with 
$M,\o,s$ replaced by $M^\sh,\o\i,s^*$); hence, by (d), it is exact (we use that $l(\o)<l(s^*\o\i)$). It 
follows that the sequence (f) is exact. From this we deduce the exactness of (g) in the same way as we have 
deduced the exactness of (c) from that of (b).

\subhead 4.5\endsubhead
Let $w\in W$. We set $N=R_w$. For $r\in R$ and $b\in N$ we write $b\cir r$ for the element of $N$ 
given by the right $R$-module structure on $N$. We define some subsets of $R.N.R$ as follows:
$$Y=\{1.\a b.1+\a.b.1+1.\a b'.\a^*+\a.b'.\a^*;b,b'\in N\},$$
$$Y'=\{1.b'\circ\a^*.1+\a.b\cir\a^*.1+1.b'.\a^*+\a.b.\a^*;b,b'\in N\},$$
$$V=\{1.\a b\cir\a^*.1+\a.b\cir\a^*.1+1.\a b.\a^*+\a.b.\a^*;b\in N\}=Y\cap Y',$$
$$Z=\{1.(\a b+b'\cir\a^*-\a b''\cir\a^*).1+\a.b.1+1.b'.\a^*+\a.b''.\a^*;b,b',b''\in N\}=Y+Y'.$$
It is easy to verify that $Y,Y'$ are subobjects of $R.N.R$ in $\car$. Hence $V,Z$ are subobjects of $R.N.R$ 
in $\car$.

By a straightforward computation we see that (a)--(d) below hold:

(a) the map $\t_1:V@>>>N[-4]$, $1.\a b\cir\a^*.1+\a.b\cir\a^*.1+1.\a b.\a^*+\a.b.\a^*\m b$ is an 
isomorphism in $\car$;

(b) the map $Y@>>>R_{ws^*}[-2]$, $1.\a b.1+\a.b.1+1.\a b'.\a^*+\a.b'.\a^*\m b-b'\cir\a^*$, induces an 
isomorphism $\t_2:Y/V@>\si>>R_{ws^*}[-2]$ in $\car$;

(c) the map $Y'@>>>R_{sw}[-2]$, $1.b'\circ\a^*.1+\a.b\cir\a^*.1+1.b'.\a^*+\a.b.\a^*\m {}^s(b'-\a b)$, 
induces an isomorphism $\t_3:Y/V'@>\si>>R_{sw}[-2]$ in $\car$;

(d) the map $R.N.R@>>>R_{sws^*}$, 

$1.b_0.1+\a.b_1.1+1.b_2.\a^*+\a.b_3.\a^*\m{}^s(-b_0+\a b_1+b_2\cir\a^*-\a b_3\cir\a^*)$,
\nl
induces an isomorphism $\t_4:R.N.R/Z@>\si>>R_{sws^*}$.

\mpb

In the remainder of this section we fix $M\in\tC$.

\proclaim{Lemma 4.6} Assume that $l(w)<l(ws^*)$. The obvious sequence 
$$0@>>>\un{\un{V^M}}@>>>\un{\un{Y^M}}@>>>\un{\un{(Y/V)^M}}@>>>0\tag a$$
is exact.
\endproclaim
We can identify $N.R[-2]=Y$ (as objects of $\car$) by $b.r\m\a.b.r+1.\a b.r$ (for $b\in N,r\in R$). We can 
identify $N[-4]=V$ via $\t_1$ in 4.5(a) and $R_{ws^*}[-2]=Y/V$ via $\t_2$ in 4.5(b). Then (a) becomes a 
sequence
$$0@>>>\un{\un{(N[-4])^M}}@>>>\un{\un{(N.R[-2])^M}}@>>>\un{\un{(R_{ws^*}[-2])^M}}@>>>0.$$
By 4.2(iii) we can identify $\un{\un{(N.R[-2])^M}}=\un{\un{(N[-4])^{M.R}}}$ (as $\RR$-vector spaces). The 
previous sequence becomes a sequence
$$0@>>>\un{\un{N^{M[4]}}}@>>>\un{\un{N^{M[4].R}}}@>>>\un{\un{(R_{ws^*})^{M[4]}[2]}}@>>>0.$$
This is of the type appearing in 4.4(c) with $M$ replaced by $M[4]$ hence is exact by 4.4(d). The lemma is 
proved.

\proclaim{Lemma 4.7} Assume that $l(sw)<l(sws^*)$. The obvious sequence
$$0@>>>\un{\un{(Z/Y)^M}}@>>>\un{\un{(R.N.R/Y)^M}}@>>>\un{\un{(R.N.R/Z)^M}}@>>>0\tag a$$
is exact.
\endproclaim
Consider the exact sequence $0@>>>N[-2]@>c>>R.N@>c'>>R_{sw}@>>>0$
in which $c$ is $b\m \a.b+\a b.1$ and $c'$ maps $r'.b$ to $r'{}^sb$ (where $r'\in R,b\in N$). Applying
$\ot_{R^{s^*}}R$ we obtain an exact sequence $0@>>>N.R[-2]@>>>R.N.R@>>>R_{sw}.R@>>>0$. Here we identify
$N.R[-2]=Y$ as in the proof of 4.6 and we obtain an exact sequence $0@>>>Y@>>>R.N.R@>>>R_{sw}.R@>>>0$ in
$\car$. Hence we obtain an identification $R.N.R/Y=R_{sw}.R$ under which $r'.b.r \in R.N.R/Y$ corresponds to
$r'{}^sb.r\in R_{sw}.R$.
We identify $Z/Y=(Y+Y')/Y=Y'/V=R_{sw}[-2]$ via the isomorphism $\t_3$ in 4.5(c) and
$R.N.R/Z=R_{sws^*}$ via the isomorphism $\t_4$ in 4.5(d). 
Then (a) becomes 
$$0@>>>\un{\un{(R_{sw}[-2])^M}}@>>>\un{\un{(R_{sw}.R)^M}}@>>>\un{\un{R_{sws^*}^M}}@>>>0$$
By 4.2(iii) we can identify $(R_{sw}.R)^M=(R_{sw}[-2])^{M.R}$. The previous sequence becomes
$$0@>>>\un{\un{(R_{sw}[-2])^M}}@>>>\un{\un{(R_{sw}[-2])^{M.R}}}@>>>\un{\un{R_{sws^*}^M}}@>>>0.$$

This sequence is (up to shift) of the type appearing in 4.4(g) (with $\o$ replaced by $sw$) hence is exact
by 4.4(h). The lemma is proved.

\proclaim{Lemma 4.8} Assume that $l(w)<l(sw)$. The obvious sequence 
$$0@>>>\un{\un{Y^M}}@>>>\un{\un{(R.N.R)^M}}@>>>\un{\un{(R.N.R/Y)^M}}@>>>0\tag a$$
is exact.
\endproclaim
We identify $Y=N.R[-2]$ as in the proof of 4.6 and $R.N.R/Y=R_{sw}.R$ as in the proof of 4.7.
Then (a) becomes the sequence
$$0@>>>\un{\un{(N.R[-2])^M}}@>>>\un{\un{(R.N.R)^M}}@>>>\un{\un{(R_{sw}.R)^M}}@>>>0.$$
By 4.2(iii), 4.3(iii) we can identify
$$(N.R[-2])^M=(N[-4])^{M.R},\qua (R.N.R)^M=(N[-4])^{R.M.R},$$
$$(R_{sw}.R)^M=(R_{sw}[-2])^{M.R}$$
and the previous sequence becomes
$$0@>>>\un{\un{(N[-4])^{M.R}}}@>>>\un{\un{(N[-4])^{R.M.R}}}@>>>\un{\un{(R_{sw}[-2])^{M.R}}}@>>>0.$$
This sequence is of the type appearing in 4.4(c) with $M$ replaced by $M.R[4]$, hence is exact by 4.4(d).
The lemma is proved.

\subhead 4.9\endsubhead
We set $P=\hom^\bul(M,N)$ regarded as an object of $\car$ as in 4.3(a).
We define subsets $\cv,\cy,\cy',\cz$ of $P$ as follows:
$$\cv=\{f\in P;f(\a m)=\a f(m),f(m\a^*)=f(m)\cir\a^*\text{ for all }m\in M\};$$
$$\cy=\{f\in P;f(\a m)=\a f(m)\text{ for all }m\in M\};$$
$$\cy'=\{f\in P;f(m\a^*)=f(m)\cir\a^*\text{ for all }m\in M\};$$
$$\cz=\{f\in P;f(\a m\a^*)-\a f(m\a^*)-f(\a m)\cir\a^*+\a f(m)\cir\a^*=0\text{ for all }m\in M\}.$$
Note that $\cv,\cy,\cy',\cz$ are subobjects of $P$ in $\car$. Under the bijection $P\lra(R.N.R)^M[4]$ in 
4.3(ii), $\cv,\cy,\cy',\cz$ correspond respectively to the subsets $V^M,Y^M,Y'{}^M,Z^M$ of $(R.N.R)^M$. Thus
we have natural bijections $\cv\lra V^M$, $\cy\lra Y^M$, $\cy'\lra Y'{}^M$, $\cz\lra Z^M$ as 
$(R^s,R^{s^*})$-bimodules. From the definitions it is clear that
$$\cv=N^M\tag a$$
as objects of $\car$. Since $P\cong N^{R.M.R}$ as objects of $\car$, we see from 2.1(a) that $P$ is a 
finitely generated right $R$-module. Since $R$ is a Noetherian ring, it follows that $\cv,\cy,\cy',\cz$ 
(which are subobjects of $P$) are also finitely generated right $R$-modules.

\proclaim{Lemma 4.10} Assume that $l(w)<l(ws^*)$. The map (in $\car$)

(a) $\cy@>>>(R_{ws^*}[2])^M$, $f\m f'$, where $f'(m)=f(m\a^*)-f(m)\cir\a^*$,
\nl
induces an isomorphism $\un{\cy/\cv}@>\si>>\un{(R_{ws^*}[2])^M}$ and an isomorphism
$(\cy/\cv)\hat{}@>\si>>((R_{ws^*}[2])^M)\hat{}$.
\endproclaim
The map (a) is clearly a well defined morphism in $\car$ and its kernel is
clearly equal to $\cv$. Thus we have an exact sequence $0@>>>\cv@>>>\cy@>>>(R_{ws^*}[2])^M$ (in $\car$).
Using 4.9(a) and the identification $\cy=N^{M.R}$ (see 4.2(i)) this exact sequence becomes an exact 
sequence $0@>>>N^M@>>>N^{M.R}@>>>(R_{ws^*}[2])^M$ (in $\car$) which induces the exact sequence 
$0@>>>\un{N^M}@>>>\un{N^{M.R}}@>>>\un{(R_{ws^*}[2])^M}@>>>0$ (a special case of 4.4(c),(d)) that is an exact
sequence $0@>>>\un{\cv}@>>>\un{\cy}@>>>\un{(R_{ws^*}[2])^M}@>>>0$. Applying $\ot_R\RR$ to the exact sequence
$0@>>>\cv@>>>\cy@>>>\cy/\cv@>>>0$ we deduce an exact sequence $\un{\cv}@>>>\un{\cy}@>>>\un{\cy/\cv}@>>>0$. 
It follows that both $\un{(R_{ws^*}[2])^M}$ and $\un{\cy/\cv}$ can be identified with the cokernel of the 
map $\un{\cv}@>>>\un{\cy}$. 
Thus, $\un{\cy/\cv}@>\si>>\un{(R_{ws^*}[2])^M}$. Now the injective homomorphism
$\cy/\cv@>>>(R_{ws^*}[2])^M$ induces an injective homomorphism
$(\cy/\cv)\hat{}@>>>((R_{ws^*}[2])^M)\hat{}$ which becomes surjective after applying $\ot_{\hR}\RR$; hence,
by the Nakayama lemma, it is surjective before applying $\ot_{\hR}\RR$. The lemma is proved.

\proclaim{Lemma 4.11} Assume that $l(w)<l(sw)$. The map (in $\car$)

(a) $\cy'@>>>(R_{sw}[2])^M$, $f\m f'$, where $f'(m)={}^s(f(\a m)-\a f(m))$,
\nl
induces an isomorphism $\un{\cy'/\cv}@>\si>>\un{(R_{sw}[2])^M}$ and an isomorphism
$(\cy'/\cv)\hat{}@>\si>>((R_{sw}[2])^M)\hat{}$.
\endproclaim
The proof is almost a repetition of that of Lemma 4.10.
The map (a) is clearly a well defined morphism in $\car$ and its kernel is
clearly equal to $\cv$. Thus we have an exact sequence $0@>>>\cv@>>>\cy'@>>>(R_{sw}[2])^M$ (in $\car$).
Using 4.9(a) and the identification $\cy'=N^{R.M}$ (see 4.2(i)) this exact sequence becomes an exact 
sequence $0@>>>N^M@>>>N^{R.M}@>>>(R_{sw}[2])^M$ (in $\car$) which induces the exact sequence 
$0@>>>\un{N^M}@>>>\un{N^{R.M}}@>>>\un{(R_{sw}[2])^M}@>>>0$ (a special case of 4.4(b),(d)) that is an exact
sequence $0@>>>\un{\cv}@>>>\un{\cy'}@>>>\un{(R_{sw}[2])^M}@>>>0$. Applying $\ot_R\RR$ to the exact sequence
$0@>>>\cv@>>>\cy'@>>>\cy'/\cv@>>>0$ we deduce an exact sequence 
$\un{\cv}@>>>\un{\cy'}@>>>\un{\cy'/\cv}@>>>0$. 
It follows that both $\un{(R_{sw}[2])^M}$ and $\un{\cy'/\cv}$ can be identified with the cokernel of the 
map $\un{\cv}@>>>\un{\cy'}$. 
Thus, $\un{\cy'/\cv}@>\si>>\un{(R_{sw}[2])^M}$. Now the injective homomorphism
$\cy'/\cv@>>>(R_{sw}[2])^M$ induces an injective homomorphism
$(\cy'/\cv)\hat{}@>>>((R_{sw}[2])^M)\hat{}$ which becomes surjective after applying $\ot_{\hR}\RR$; hence,
by the Nakayama lemma, it is surjective before applying $\ot_{\hR}\RR$. The lemma is proved.

\proclaim{Lemma 4.12} Assume that $l(w)<l(sw)$. Let $P'=\hom'{}^\bul(M,R_{sw})$; we view $P'$ as an object 
of $\car$ as in 4.2(i). The map (in $\car$)

(a) $P@>>>P'[2]$, $f\m f'$, $f'(m)={}^s(f(\a m)-\a f(m))$,
\nl
induces an isomorphism $\un{P/\cy}@>\si>>\un{P'[2]}$ (hence, using $P'=R_{sw}^{M.R}$, see 4.2(i)) an 
isomorphism $\un{P/\cy}@>\si>>\un{(R_{sw}[2])^{M.R}}$; it also induces an isomorphism
$(P/\cy)\hat{}@>\si>>((R_{sw}[2])^{M.R})\hat{}$.
\endproclaim
The map (a) is clearly a well defined morphism in $\car$ and its kernel is clearly equal to $\cy$. Thus we 
have an exact sequence $0@>>>\cy@>>>P@>>>P'[2]$ in $\car$. By 4.2(i) we can identify $\cy=N^{M.R}$ and our 
exact sequence becomes the exact sequence 
$0@>>>N^{M.R}@>>>N^{R.M.R}@>>>(R_{sw}[2])^{M.R}$ in $\car$ which induces an exact sequence
$0@>>>\un{N^{M.R}}@>>>\un{N^{R.M.R}}@>>>\un{(R_{sw}[2])^{M.R}}@>>>0$ 
(a special case of 4.4(b),(d) with $M$ replaced by $M.R$). Thus we have an exact sequence
$0@>>>\un{\cy}@>>>\un{P}@>>>\un{P'[2]}@>>>0$.
Applying $\ot_R\hR$ to the exact sequences 

$0@>>>\cy@>>>P@>>>P/\cy@>>>0$, $0@>>>\cy@>>>P@>>>P'[2]$,
\nl
 we obtain exact sequences 

$0@>>>\hat\cy@>>>\hat P@>>>\widehat{P/\cy}@>>>0$, $0@>>>\hat\cy@>>>\hat P@>>>\hat P'[2]$. 
\nl
From the surjectivity of $\un{P}@>>>\un{P'[2]}$ and the Nakayama 
lemma it follows that $\hat P@>>>\hat P'[2]$ in the last exact sequence is surjective.
Hence the obvious map $\widehat{P/\cy}@>>>\hat P'[2]$ is an isomorphism (both sides can be identified with
$\coker(\hat\cy@>>>\hat P)$). Applying $\ot_{\hR}\RR$ we deduce that the obvious map
$\un{P/\cy}@>>>\un{P'[2]}$ is an isomorphism. 
Thus, $\un{P/\cy}@>\si>>\un{(R_{sw}[2])^{M.R}}$.
Now the injective homomorphism
$P/\cy@>>>(R_{sw}[2])^M$ induces an injective homomorphism
$(P/\cy)\hat{}@>>>((R_{sw}[2])^M)\hat{}$ which becomes surjective after applying $\ot_{\hR}\RR$; hence,
by the Nakayama lemma, it is surjective before applying $\ot_{\hR}\RR$. The lemma is proved.

\proclaim{Lemma 4.13} Assume that $l(w)<l(sw)<l(sws^*)$. The map (in $\car$) $P@>>>(R_{sws^*}[4])^M$, 

(a) $f\m f'$, $f'(m)={}^s(f(\a m\a^*)-\a f(m\a^*)-f(\a m)\cir\a^*+\a f(m)\cir\a^*)$, 
\nl
induces an isomorphism $\un{P/\cz}@>\si>>\un{(R_{sws^*}[4])^M}$ and an isomorphism
$(P/\cz)\hat{}@>\si>>((R_{sws^*}[4])^M)\hat{}$.

The map (in $\car$) $\cz@>>>(R_{sw}[2])^M$, 

(b) $f\m f'$, $f'(m)={}^s(f(\a m)-\a f(m))$
\nl
induces an isomorphism $\un{\cz/\cy}@>\si>>\un{(R_{sw}[2])^M}$ and an isomorphism
$(\cz/\cy)\hat{}@>\si>>((R_{sw}[2])^M)\hat{}$.
\endproclaim
Th map (a) is clearly a well defined morphism in $\car$ and its kernel is
clearly equal to $\cz$. Thus we have an exact sequence $0@>>>\cz/\cy@>>>P/\cy@>>>(R_{sws^*}[4])^M$.
Applying $\ot_R\hR$ gives again an exact sequence
$$0@>>>\widehat{\cz/\cy}@>>>\widehat{P/\cy}@>>>((R_{sws^*}[4])^M)\hat{}.\tag c$$
From 4.4(f),(h) we have an exact sequence
$$0@>>>\un{(R_{sw}[2])^M}@>>>\un{(R_{sw}[2])^{M.R}}@>>>\un{(R_{sws^*}[4])^M}@>>>0.\tag d$$
Hence $\un{(R_{sw}[2])^{M.R}}@>>>\un{(R_{sws^*}[4])^M}$ is surjective, that is (using 4.12)
$\un{P/\cy}@>>>\un{R_{sws^*}^M}$ is surjective.
Using this and Nakayama lemma we see that $\widehat{P/\cy}@>>>(R_{sws^*}^M)\hat{}$ is surjective.
This is just the last map in (c); thus, (c) becomes an exact sequence
$$0@>>>\widehat{\cz/\cy}@>>>\widehat{P/\cy}@>>>((R_{sws^*}[4])^M)\hat{}@>>>0.$$
This exact sequence of $\hR$-modules splits since, by 2.1(a), the $\hR$-module\lb
$((R_{sws^*}[4])^M)\hat{}$ is free. Hence, applying $\ot_{\hR}\RR$ gives an exact sequence
$$0@>>>\un{\cz/\cy}@>>>\un{P/\cy}@>>>\un{(R_{sws^*}[4])^M}@>>>0.\tag e$$
From the obvious exact sequence $0@>>>\cz/\cy@>>>P/\cy@>>>P/\cz@>>>0$ we deduce an exact sequence
$\un{\cz/\cy}@>>>\un{P/\cy}@>>>\un{P/\cz}@>>>0$.
Using this and (d) we see that both $\un{P/\cz}$ and $\un{(R_{sws^*}[4])^M}$ can be identified with the
cokernel of the map $\un{\cz/\cy}@>>>\un{P/\cy}$. 
Using (d) and (e), where we identify
$\un{(R_{sw}[2])^{M.R}}=\un{P/\cy}$ (see 4.12), we see that both $\un{\cz/\cy}$ and $\un{(R_{sw}[2])^M}$ can 
be identified with the kernel of the map $\un{P/\cy}@>>>\un{P/\cz}$. 
Thus, we have $\un{P/\cz}@>\si>>\un{(R_{sws^*}[4])^M}$ and $\un{\cz/\cy}@>\si>>\un{(R_{sw}[2])^M}$.
Now, the injective homomorphism $P/\cz@>>>(R_{sws^*}[4])^M$ (resp. $\cz/\cy@>>>(R_{sw}[2])^M$)
induces an injective homomorphism
$(P/\cz)\hat{}@>>>((R_{sws^*}[4])^M)\hat{}$ (resp. $(\cz/\cy)\hat{}@>>>((R_{sw}[2])^M)\hat{}$)  
which becomes surjective after applying $\ot_{\hR}\RR$; hence,
by the Nakayama lemma, it is surjective before applying $\ot_{\hR}\RR$. The lemma is proved.

\proclaim{Lemma 4.14} Assume that $l(w)<l(sw)=l(ws^*)$. The obvious sequence
$$0@>>>\un{\cv}@>>>\un{P}@>>>\un{P/\cv}@>>>0$$
is exact.
\endproclaim
From the exact sequence $0@>>>\cy/\cv@>>>P/\cv@>>>P/\cy@>>>0$ we deduce an exact sequence
$0@>>>(\cy/\cv)\hat{}@>>>(P/\cv)\hat{}@>>>(P/\cy)\hat{}@>>>0$ in which 
$(\cy/\cv)\hat{}$ is a free $\hR$-module (by 4.10 and 2.1(a)) and
$(P/\cy)\hat{}$ is a free $\hR$-module (by 4.12 and 2.1(a)). It follows that

(a) $(P/\cv)\hat{}$ is a free $\hR$-module.
\nl
From the obvious exact sequence $0@>>>\cv@>>>P@>>>P/\cv@>>>0$ we deduce an exact sequence
$0@>>>\hat\cv@>>>\hat P@>>>(P/\cv)\hat{}@>>>0$ which is split, due to (a). It follows that it remains exact
after applying $\ot_{\hR}\RR$. The lemma is proved.

\proclaim{Lemma 4.15} Assume that $l(w)<l(sw)<l(sws^*)$. The obvious sequence 
$$0@>>>\un{\cz/\cv}@>>>\un{P/\cv}@>>>\un{P/\cz}@>>>0$$
is exact.
\endproclaim
From the obvious exact sequence $0@>>>\cz/\cv@>>>P/\cv@>>>P/\cz@>>>0$ we deduce an exact sequence
$0@>>>(\cz/\cv)\hat{}@>>>(P/\cv)\hat{}@>>>(P/\cz)\hat{}@>>>0$ which is split, since the $\hR$-module 
$(P/\cz)\hat{}$ is free, by 4.13 and 2.1(a). It follows that it remains exact after applying $\ot_{\hR}\RR$. 
The lemma is proved.

\proclaim{Lemma 4.16} Assume that $l(w)<l(sw)<l(sws^*)$. The sum of the obvious homomorphisms
$\un{\cy/\cv}@>c>>\un{\cz/\cv}$ and $\un{\cy'/\cv}@>c'>>\un{\cz/\cv}$ is an isomorphism 
$\un{\cy/\cv}\op\un{\cy'/\cv}@>\si>>\un{\cz/\cv}$.
\endproclaim
From the obvious exact sequence 
$0@>>>\cy/\cv@>>>\cz/\cv@>>>\cz/\cy@>>>0$ we deduce an exact sequence
$0@>>>(\cy/\cv)\hat{}@>>>(\cz/\cv)\hat{}@>>>(\cz/\cy)\hat{}@>>>0$ which is split, since the $\hR$-module 
$(\cz/\cy)\hat{}$ is free, by 4.13 and 2.1(a). It follows that after applying $\ot_{\hR}\RR$ we get an
exact sequence 
$$0@>>>\un{\cy/\cv}@>c>>\un{\cz/\cv}@>d>>\un{\cz/\cy}@>>>0.$$
We consider the composition $dc':\un{\cy'/\cv}@>>>\un{\cz/\cy}$. By 4.11 we can identify
$\un{\cy'/\cv}=\un{(R_{sw}[2])^M}$ and by 4.13 we can identify $\un{\cz/\cy}=\un{(R_{sw}[2])^M}$.
Under these identifications the map $dc'$ becomes the identity map of $\un{(R_{sw}[2])^M}$.
In particular, $dc'$ is an isomorphism. This implies immediately the statement of the lemma.

\head 5. Trace computations\endhead
\subhead 5.1\endsubhead
To simplify notation, for $x\in W,r\in R$ we shall write ${}^xr^*$ instead of ${}^x(r^*)$.
Recall that if $x\in\II$, then $r\m{}^xr^*$ is an involution $R@>>>R$ denoted by $\ff_x$ in 3.2.

In this section we fix $(M,\ph)\in\tC_\sh$, $s\in\II$ and $w\in\II$ such that $l(w)<l(sw)$; we have 
automatically $l(w)<l(ws^*)$. As in 4.5 we set $N=R_w$. The notation $b\cir r$ for $b\in N,r\in R$ is as in 
4.5. In the case where $sw=ws^*$ we set $N'=R_{sw}$. In the case where $sw\ne ws^*$ we set $N''=R_{sws^*}$.

For $b\in N,r\in R$ we have 
$$\ff_w(b\cir r)=r^*\ff_w(b),\qua \ff_w(rb)=\ff_w(b)\cir r^*.$$
The involution $f\m f^!$, $N^M@>>>N^M$, given by $f^!(m)=\ff_w(f(\ph(m))$ 
induces an involution $\Th:\un{N^M}@>>>\un{N^M}$.
In the case where $sw=ws^*$, we have $sw\in\II$ and the involution $f\m f^!$, $N'{}^M@>>>N'{}^M$, given by 
$f^!(m)=\ff_{sw}(f(\ph(m))$ induces an involution $\Th':\un{N'{}^M}@>>>\un{N'{}^M}$.
In the case where $sw\ne ws^*$, we have $sws^*\in\II$ and the involution $f\m f^!$, $N''{}^M@>>>N''{}^M$, 
given by $f^!(m)=\ff_{sws^*}(f(\ph(m))$ induces an involution $\Th'':\un{N''{}^M}@>>>\un{N''{}^M}$.
Now $\Th$ (or $\Th'$ or $\Th''$, if defined) induces a degree preserving involution of $\un{N^M}$ (or 
$\un{N'{}^M}$, or $\un{N''{}^M}$) denoted again by $\Th$ (or $\Th'$ or $\Th''$).

By 3.6 we have $(R.M.R,\ph')\in\tC_\sh$ where $\ph':R.M.R@>>>R.M.R$ is the $\RR$-linear map such that 
$r_1.m.r_2\m r_2^*.\ph(m').r_1^*$ for $r_1,r_2\in R$, $m\in M$. (Recall that $R.M.R\in\tC$ is defined in 
4.3.) We have $\ph'{}^2=1$. Let $\Ps:N^{R.M.R}@>>>N^{R.M.R}$ be the $\RR$-linear involution such that
for any $F\in N^{R.M.R}$ and any $r_1,r_2\in R$, $m\in M$, we have 
$$\Ps(F)(r_1.m.r_2)=\ff_w(F(\ph'(r_1.m.r_2)))=\ff_w(F(r_2^*.\ph(m).r_1^*)).$$
(This is a special case of the definition of $f\m f^!$ in 1.1.). 
It induces a degree preserving involution of $\un{N^{R.M.R}}$ denoted again by $\Ps$. 

We now state the main result of this section. (In this section all traces are taken over $\RR$.)

\proclaim{Theorem 5.2} Recall that $w\in\II$, $l(w)<l(sw)$. Let $i\in\ZZ$. If $sw\ne ws^*$ then
$$\tr(\Ps,\un{N^{R.M.R}}_i)=\tr(\Th,\un{N^M}_i)+\tr(\Th'',\un{N''{}^M}_{i+4}).\tag a$$
If $sw=ws^*$ then
$$\align&\tr(\Ps,\un{N^{R.M.R}}_i)=\tr(\Th,\un{N^M}_i)+\tr(\Th,\un{N^M}_{i+2})\\&
-\tr(\Th',\un{N'{}^M}_{i+2})+\tr(\Th',\un{N'{}^M}_{i+4}).\tag b\endalign$$
\endproclaim
Note that the following identities (with $\ph'$ as in 5.1) are equivalent to the theorem.
$$\e^w(R.M.R,\ph')=\e^w(M,\ph)+\e^{sws^*}(M,\ph)v^4 \text{ if }sw\ne ws^*,\tag c$$
$$\e^w(R.M.R,\ph')=\e^w(M,\ph)(v^2+1)+\e^{sw}(M,\ph)(v^4-v^2) \text{ if }sw=ws^*.\tag d$$
The proof will occupy the remainder of this section.

\subhead 5.3\endsubhead
We identify $N^{R.M.R}$ with $P=\hom^\bul(M,N)$ (as objects of $\car$) as in 4.3(i). Then $\Ps$ becomes an 
involution of $P$ denoted again by $\Ps$. It is given by $f\m f^!$ where $f^!(m)=\ff_w(f(\ph(m)))$. This 
induces a degree preserving involution of $\un{P}$ denoted again by $\Ps$. For any $i$ we have clearly 
$$\tr(\Ps,\un{N^{R.M.R}}_i)=\tr(\Ps,\un{P}_i).\tag a$$ 

\subhead 5.4\endsubhead
In this subsection we assume that $sw\ne ws^*$ so that $l(w)<l(sw)<l(sws^*)$ and $sws^*\in\II$.
Let $\cv,\cy,\cy',\cz$ be the subobjects of $P$
defined in 4.9. From the definition we see that $\Ps:P@>>>P$ preserves $\cv$ and $\cz$; it interchanges
$\cy$ and $\cy'$. Now for any $\x\in P$ we have $\Ps(\x R^{>0})=R^{>0}\Ps(\x)=\Ps(\x)R^{>0}$. 
(We use that $R^{>0}b=b\cir R^{>0}$ for $b\in N$.) It follows that $\cv R^{>0},\cz R^{>0}$ are preserved
by $\Ps$ and $\cy R^{>0},\cy'R^{>0}$ are interchanged by $\Ps$. Hence $\Ps$ induces involutions of 
$\un{\cv}_i,\un{P/\cz}_i,\un{\cz/\cv}_i$ (denoted again by $\Ps$) and the two summands $\un{\cy/\cv}_i$,
$\un{\cy'/\cv}_i$ of $\un{\cz/\cv}_i$ (see 4.16) are intechanged by $\Ps:\un{\cz/\cv}_i@>>>\un{\cz/\cv}_i$.
Hence we have $\tr(\Ps,\un{\cz/\cv}_i)=0$ and (using 4.14, 4.15) we have
$$\tr(\Ps,\un{P}_i)=\tr(\Ps,\un{\cv}_i)+\tr(\Ps,\un{P/\cv}_i)=\tr(\Ps,\un{\cv}_i)+\tr(\Ps,\un{P/\cz}_i).
\tag a$$
We now show that the map (say $\t$), $P@>>>R_{sws^*}[4]$ in 4.13(a) satisfies
$$\t(\Ps(f))=\Th''(\t(f))\tag b$$
for any $f\in P$. For $m\in M$ we have   
$$\align&\t(\Ps(f))(m)=\\&{}^s(\ff_w(f(\ph(\a m\a^*)))
-\a \ff_w(f(\ph(m\a^*)))-\ff_w(f(\ph(\a m)))\cir\a^*+\a \ff_w(f(\ph(m)))\cir\a^*),\endalign$$
$$\align&\Th''(\t(f))(m)=\ff_{sws^*}(\t(f)(\ph(m)))\\&=
\ff_{sws^*}({}^s(f(\a \ph(m)\a^*)-\a f(\ph(m)\a^*)-f(\a\ph(m))\cir\a^*+\a f(\ph(m))\cir\a^*)).\endalign$$
It is enough to show that for any $m'\in M$ we have
$${}^s(\ff_w(f(m')))=\ff_{sws^*}({}^s(f(m')))$$
or that
$${}^s({}^w(f(m')^*))={}^{sws^*}({}^{s^*}(f(m')^*)).$$
This is clear; (b) is proved. Using (b) and 4.13 we deduce that
$$\tr(\Ps,\un{P/\cz}_i)=\tr(\Th'',\un{N''{}^M}_{i+4}).\tag c$$
We have clearly $\cv=N^M$ and $\tr(\Ps,\un{\cv}_i)=\tr(\Th,\un{N^M}_i)$.
Introducing this and (c) into (a) and using 5.3(a) we obtain 5.2(a) and (equivalently) 5.2(c).

\subhead 5.5\endsubhead
In the remainder of this section we assume that $sw=ws^*$ so that $sw\in\II$. 
Note that we have ${}^w\a^*=\a$, hence $b\cir\a^*=\a b$ for $b\in N$.

In this case the involution $\Ps:P@>>>P$ preserves $PR^{s^*,>0}$ 
(more precisely, $\Ps(fR^{s^*,>0})=\Ps(f)R^{s^*,>0}$ for any $f\in P$) hence it induces an involution of 
$\un{\un{P}}$ denoted again by $\Ps$. (We use that ${}^w(R^{s^*})=R^s$ hence
${}^w(R^{s^*}\cap R^{>0})=R^{s}\cap R^{>0}$). Moreover the involution $\Ps$ of $\un{\un{P}}$ is 
$\un{\un{R}}$-linear. (We use that ${}^w\a^*=\a$.)

Let $\Ph:(R.N.R)^M@>>>(R.N.R)^M$ be the $\RR$-linear involution which corresponds to $\Ps:P@>>>P$ under the 
bijection $P[-4]@>\si>>(R.N.R)^M$ in 4.3(ii). Since that bijection is compatible with the 
$(R^s,R^{s^*})$-bimodule structures, it follows that $\Ph$ preserves the subset $(R.N.R)^MR^{s^*,>0}$, more
precisely we have

(a) $\Ph(\x R^{s^*,>0})=\Ph(\x)R^{s^*,>0}$ for any $\x\in(R.N.R)^M$,
\nl
hence $\Ph$ induces an $\RR$-linear involution of $\un{\un{(R.N.R)^M}}$ (which is not necessarily
$\un{\un{R}}$-linear).
For any $i$ we have from the definition:
$$\tr(\Ph,\un{\un{(R.N.R)^M}}_i)=\tr(\Ps,\un{\un{P}}_{i-4}).\tag b$$
Note that $\un{\un{P}}$ is a free right $\RR[\a^*]/(\a^{*2})$-modules.
Hence we have exact an sequence of $\RR$-vector spaces
$$0@>>>\un{P}_{i-6}@>c>>\un{\un{P}}_{i-4}@>d>>\un{P}_{i-4}@>>>0,$$
where $c$ is induced by right multiplication by $\a^*$ and we have $d\Th=\Th d$, $c\Th=\Th c$. It follows 
that we have
$$\tr(\Ps,\un{\un{P}}_{i-4})=\tr(\Ps,\un{P}_{i-4})+\tr(\Th,\un{P}_{i-6}).$$
Introducing this in (b) we obtain
$$\tr(\Ph,\un{\un{(R.N.R)^M}}_i)=\tr(\Ps,\un{P}_{i-4})+\tr(\Ps,\un{P}_{i-6}).\tag c$$

\subhead 5.6\endsubhead
We define a map $\x\m\che\x$, $R.N.R@>>>R.N.R$ by 
$$\align&1.b_0.1+\a.b_1.1+1.b_2.\a^*+\a.b_3.\a^*\\&
\m 1.\ff_w(b_0).1+\a.\ff_w(b_2).1+1.{}^w\ff_w(b_1).\a^*+\a.\ff_w(b_3).\a^*\tag a\endalign$$
where $b_i\in N$. Then $\x\m\che\x$ is an involution of the $\RR$-vector space $R.N.R$
such that $(r_1.b.r_2)\che{}=r_2^*.\ff_w(b).r_1^*$ for $r_1,r_2\in R,b\in N$. Hence in the $(R,R)$-bimodule 
structure of $R.N.R$ we have $(r\x)\che{}=\che\x r^*$, $(\x r)\che{}=r^*\che\x$ for $r\in R$, 
$\x\in R.N.R$. Thus $(R.N.R,\x\m\che\x)\in\car_\sh$.

From the definitions we see that $\Ph:(R.N.R)^M@>>>(R.N.R)^M$ is given explicitly by $G\m G^!$ where for any
$G\in (R.N.R)^M$ and any $m\in M$ we have
$$G^!(m)=(G(\ph(m)))\che{}.\tag b$$
(This is a special case of the definition of $f\m f^!$ in 1.1.)

\subhead 5.7\endsubhead
Let $V,Y$ be the subsets of $R.N.R$ defined as in 4.5. (They are subobjects in $\car$.)
In addition to the subsets $V,Y$ we shall need the following subsets of $R.N.R$:
$$X=\{1.\a b'\cir\a^*.1+\a.b.1+1.b.\a^*+\a.b'.\a^*;b,b'\in N\},$$
$$U=\{1.(\a b-b'\cir\a^*+\a b''\cir\a^*).1+\a.b.1+1.b'.\a^*+\a.b''.\a^*;b,b',b''\in N\}=X+Y.$$
Note that $X\cap Y=V$. Using our assumptions on $w$, it is easy to verify that $X$ is a subobject of $R.N.R$
in $\car$ hence $U=X+Y$ is a subobject of $R.N.R$ in $\car$.

By a straightforward computation we see that (a),(b) below hold:

(a) the map $X@>>>N[-2]$, $1.\a b'\cir\a^*.1+\a.b.1+1.b.\a^*+\a.b'.\a^*\m{}^s(b-\a b')$ induces an 
isomorphism $X/V@>\si>>N[-2]$ in $\car$;

(b) the map $R.N.R@>>>N'$, 

$1.b_0.1+\a.b_1.1+1.b_2.\a^*+\a.b_3.\a^*\m{}^s(-b_0+\a b_1-\a b_2+\a^2b_3)$,
\nl
induces an isomorphism $R.N.R/U@>\si>>N'$ in $\car$.

\proclaim{Lemma 5.8} The obvious sequence 
$$0@>>>\un{\un{(U/Y)^M}}@>>>\un{\un{(R.N.R/Y)^M}}@>>>\un{\un{(R.N.R/U)^M}}@>>>0\tag a$$
is exact.
\endproclaim
Note that $N'.R=N.R$. Indeed, it is enough to show that $N=N'$ as $(R,R^{s^*})$-bimodules. It is also enough
to show that if $r\in R^{s^*}$, then ${}^wr={}^{ws^*}r$; this follows from ${}^{s^*}r=r$.
We identify $R.N.R/Y=R_{sw}.R=N'.R$ (hence $R.N.R/Y=N.R$) as in the proof of 4.7, $U/Y=(X+Y)/Y=X/V=N[-2]$ as 
in 5.1(a) and $R.N.R/U=N'$ as in 5.1(b). Then (a) becomes a sequence
$$0@>>>\un{\un{(N[-2])^M}}@>>>\un{\un{(N.R)^M}}@>>>\un{\un{N'{}^M}}@>>>0$$
or equivalently a sequence
$$0@>>>\un{\un{(N[-2])^M}}@>>>\un{\un{(N[-2])^{M.R}}}@>>>\un{\un{N'{}^M}}@>>>0$$
which is exact by 4.4(g),(h).
The lemma is proved.

\subhead 5.9\endsubhead 
We write $N^0=R^s$, $N^1=\a R^s=R^s\cir\a^*$ so that $N=N^0\op N^1$ as an $(R^s,R^{s^*})$-bimodule.
For $b\in N$ we can write uniquely $b=b^0+b^1$ where $b^i\in N^i$. For $i\in\{0,1\}$ let 
$$(R.N.R)^i=\{1.b_0.1+\a.b_1.1+1.b_2.\a^*+\a.b_3.\a^*\in R.N.R;b_i\in N^i\text{ for }i=0,1,2,3\}.$$
Using the fact that $\a^2\in R^s$ we see that $(R.N.R)^i$ is a subobject of $R.N.R$ in $\car$. Thus, we have
$R.N.R=(R.N.R)^0\op(R.N.R)^1$ as objects of $\car$. For $i\in\{0,1\}$ we set
$$X^i=X\cap(R.N.R)^i=\{1.\a b'\cir\a^*.1+\a.b.1+1.b.\a^*+\a.b'.\a^*;b,b'\in N^i\},$$
$$(R.N.R/X)^i=(R.N.R)^i/X^i.$$
Then $X^i$ is a subobject of $X$ in $\car$, $(R.N.R/X)^i$ is a subobject of $R.N.R/X$ in $\car$ and we have
$$X=X^0\op X^1,\qua R.N.R/X=(R.N.R/X)^0\op(R.N.R/X)^1\tag a$$
as objects of $\car$. We have
$$X=V\op X^0,\tag b$$
$$R.N.R/X=U/X\op(R.N.R/X)^0.\tag c$$
We prove (b). We must show that for any $b,b'\in N$ there are unique $\b\in N^0,\b'\in N^0, b''\in N$ such 
that
$$\align&1.\a b'\cir\a^*.1+\a.b.1+1.b.\a^*+\a.b'.\a^*=1.\a\b'\cir\a^*.1+\a.\b.1+1.\b.\a^*+\a.\b'.\a^*\\&+
1.\a b''\cir\a^*.1+\a.\a b''.1+\a.b''.\a^*+1.\a b''.\a^*\endalign$$
or equivalently $b=\b+\a b'',b'=\b'+b''$. Setting $b''=b-\b'$ we see that we must show that there are unique
$\b\in N^0,\b'\in N^0$ such that $b-\a b'=\b-\a\b'$. This is obvious.

We prove (c). It is enough to show that 

(i) $R.N.R=U+(R.N.R)^0$,

(ii) $U\cap((R.N.R)^0+X^1)=X$.
\nl
For (i) we must show that given $b_1,b_2,b_3,b_4\in N$ there exist $b,b',b''\in N$ and 
$\b_1,\b_2,\b_3,\b_4\in N^0$ such that
$$b_1=b+\b_1, b_2=b'+\b_2,b_3=b''+\b_3, b_4=\a b-b'\cir\a^*+\a b''\cir\a^* +\b_4.$$
Setting $\b_2=\b_3=0$, $b=b_1-\b_1,b'=b_2, b''=b_3$ we see that it is enough to show that there exist 
$\b_1,\b_4\in N^0$ such that
$$b_4-\a b_1+\a b_2-\a^2 b_3=\b_4-\a\b_1.$$
This is obvious.

For (ii) we must show that given $b,b',b''\in N$ and $\b,\b'\in N^1$ such that
$$\align&1.(\a b-b'\cir\a^*+\a b''\cir\a^*).1+\a.b.1+1.b'.\a^*+\a.b''.\a^*\\&
-(1.\a\b'\cir\a^*.1+\a.\b.1+1.\b.\a^*+\a.\b'.\a^*)\in (R.N.R)^0,\endalign$$
we have $b=b'$. Our assumption implies $b^1=\b$, $b'{}^1=\b$, $b''{}^1=\b'$,
$(\a b-\a b'+\a^2 b'')^1=\a^2\b'$ (that is $(b-b'+\a b'')^0=\a\b'$). Thus, $(b-b')^1=0$ and $(b-b')^0=0$, so 
that $b-b'=0$. This proves (c).

Now (a),(b) yield isomorphisms (in $\car$) $X^0@>>>X/V$, $X^1@>\si>>V$; the first one is induced by the
identity map, the second one is the restriction to $X^1$ of the first projection $X=V\op X^0@>>>V$.
Moreover, (a),(c) yield isomorphisms (in $\car$) $(R.N.R/X)^0@>>>R.N.R/U$, $(R.N.R/X)^1@>>>U/X$; the first 
one is induced by the identity map, the second one is the restriction to $(R.N.R)^1$ of the first projection
$R.N.R/X=U/X\op(R.N.R/X)^0@>>>U/X$.

\proclaim{Lemma 5.10} The obvious sequence
$$0@>>>\un{\un{X^M}}@>>>\un{\un{(R.N.R)^M}}@>>>\un{\un{(R.N.R/X)^M}}@>>>0\tag a$$
is exact.
\endproclaim
Consider the obvious commutative diagram with exact horizontal and vertical lines
$$\CD
{} @.   0    @.   0   @.  0    {}\\
@.    @VVV     @VVV   @VVV   @.    \\  
0@>>> V@>>>X@>>>U/Y@>>>0 \\
@.    @VVV     @VVV   @VVV   @.    \\  
0@>>>Y @>>>R.N.R@>>>R.N.R/Y@>>>0\\
@.    @VVV     @VVV   @VVV   @.   \\   
0@>>>Y/V@>>>R.N.R/X@>>>R.N.R/U@>>>0\\
@.    @VVV     @VVV   @VVV   @.      \\
{} @.   0    @.   0   @.  0    {}\\ \endCD$$
Here the non-middle horizontal maps are split as exact sequences in $\car$. Indeed by the results in 5.9
they can be identified with the obvious split exact sequences  
$$0@>>> X^1@>>>X^0\op X^1@>>>X^0@>>>0,$$
$$0@>>>(R.N.R)^1@>>>(R.N.R)^0\op(R.N.R)^1@>>>(R.N.R)^0@>>>0.$$
From this we deduce the commutative diagram
$$\CD  
{} @.   0    @.   0   @.  0    {}\\
@.    @VVV     @VVV   @VVV   @.    \\  
0@>>> \un{\un{V^M}}@>>>\un{\un{X^M}}@>>>\un{\un{(U/Y)^M}}@>>>0 \\
@.    @VVV     @VVV   @VVV   @.    \\  
0@>>>\un{\un{Y^M}} @>>>\un{\un{(R.N.R)^M}}@>>>\un{\un{(R.N.R/Y)^M}}@>>>0\\
@.    @VVV     @VVV   @VVV   @.   \\   
0@>>>\un{\un{(Y/V)^M}}@>>>\un{\un{(R.N.R/X)^M}}@>>>\un{\un{(R.N.R/U)^M}}@>>>0\\
@.    @VVV     @VVV   @VVV   @.      \\
{} @.   0    @.   0   @.  0    {}\\  \endCD$$
in which the middle horizontal line and the non-middle vertical lines are exact sequences (see 4.8, 4.6, 
5.8) and in which the non-middle horizontal lines are (split) exact sequences. This implies, by diagram 
chasing, that the middle vertical line is an exact sequence. The lemma is proved.

\subhead 5.11\endsubhead
From 5.6(a) we see that $X$, $(R.N.R)^0$, $(R.N.R)^1$ are stable under 
the involution $\x\m\che\x$ of $R.N.R$. Hence that involution induces involutions on
$X$, on $(R.N.R)/X$, on $X^i$ and on $(R.N.R/X)^i$ (for $i=0,1$) which in turn induce (by formulas like 
5.6(b)) involutions on $X^M$, $((R.N.R)/X)^M$, on $(X^i)^M$ and on $((R.N.R/X)^i)^M$ (for $i=0,1$) which are 
denoted again by $\Ph$. Using 5.5(a) we see that each of these involutions preserve the image of
right multiplication by $R^{s^*,>0}$ hence we have induced involutions on
$\un{\un{X^M}}$, $\un{\un{((R.N.R)/X)^M}}$, on $\un{\un{(X^i)^M}}$ and on $\un{\un{((R.N.R/X)^i)^M}}$ (for 
$i=0,1$) which are denoted again by $\Ph$. 

Using the definitions we see that the exact sequence 5.10(a) is compatible with the involutions $\Ph$ on 
each of its terms. Using the definitions we also see that the obvious direct sum decompositions
$$\un{\un{X^M}}=\un{\un{(X^0)^M}}\op\un{\un{(X^1)^M}},$$
$$\un{\un{(R.N.R/X)^M}}=\un{\un{((R.N.R/X)^0)^M}}\op\un{\un{((R.N.R/X)^1)^M}}$$
are compatible with the involutions $\Ph$ on each of their terms. It follows that for $i\in\ZZ$ we have
$$\align&\tr(\Ph,\un{\un{(R.N.R)^M}}_i)=\tr(\Ph,\un{\un{X^M}}_i)+\tr(\Ph,{\un{(R.N.R/X)^M}}_i)\\&
=\tr(\Ph,\un{\un{(X^0)^M}}_i)+\tr(\Ph,\un{\un{(X^1)^M}}_i)\\&+
\tr(\Ph,\un{\un{((R.N.R/X)^0)^M}}_i)+\tr(\Ph,\un{\un{((R.N.R/X)^1)^M}}_i).\tag a\endalign$$ 

\subhead 5.12\endsubhead
By a straightforward computation we see that the maps $t_i$ in (a)--(d) below are isomorphisms in $\car$:

(a) $t_1:X^1@>>>N[-4]$, $1.\a\b'\cir\a^*.1+\a.\b.1+1.\b.\a^*+\a.\b'.\a^*\m\a\i\b+\b'$ with $\b,\b'\in N^1$;

(b) $t_1:X^0@>>>N[-2]$, $1.\a\b'\cir\a^*.1+\a.\b.1+1.\b.\a^*+\a.\b'.\a^*\m\b+\a\b'$ with $\b,\b'\in N^0$;

(c) $t_3:(R.N.R/X)^1@>>>N'[-2]$ induced by

$1.\b_0.1+\a.\b_1.1+1.\b_2.\a^*+\a.\b_3.\a^*\m\b_1-\b_2+\a\i\b_0-\a\b_3$
\nl
with $\b_0,\b_1,\b_2,\b_3\in N^1$;

(d) $t_4:(R.N.R/X)^0@>>>N'$ induced by

$1.\b_0.1+\a.\b_1.1+1.\b_2.\a^*+\a.\b_3.\a^*\m\a\b_1-\a\b_2+\b_0-\a^2\b_3$
\nl
with $\b_0,\b_1,\b_2,\b_3\in N^0$.

\subhead 5.13\endsubhead
The identities (a)--(d) below express a connection between $\x\m\che\x$ and the isomorphisms $t_j$ in 5.12:

(a) if $\x\in X^1$ then $t_1(\che\x)={}^w(t_1(\x))^*$;

(b) if $\x\in X^0$ then $t_2(\che\x)={}^w(t_2(\x))^*$;

(c) if $\x\in(R.N.R/X)^1$ then $t_3(\che\x)={}^{sw}(t_3(\x))^*$;

(d) if $\x\in(R.N.R/X)^0$ then $t_4(\che\x)={}^{sw}(t_4(\x))^*$.
\nl
Here  $t_i(\x)$ is viewed as an element of $R$ and the shift is ignored.

We prove (a). Let $\x=1.\a\b'\cir\a^*.1+\a.\b.1+1.\b.\a^*+\a.\b'.\a^*\in X^1$ be as in 5.3(a). Then 
$$\che\x=1.\a{}^w\b'{}^*\cir\a^*.1+\a.{}^w\b^*.1+1.{}^w\b^*.\a^*+\a.{}^w\b'{}^*.\a^*$$
and we must show that $\a\i{}^w\b^*+{}^wb'{}^*={}^w(\a\i\b+\b')^*$. This follows from the equality
${}^w\a^*=\a$. 

We prove (b). In this case we must show that ${}^w\b^*+\a{}^wb'{}^*={}^w(\b+\a\b')^*$ for $\b,\b'\in N^0$. 
This again follows from the equality ${}^w\a^*=\a$. 

We prove (c). Let $\x=1.\b_0.1+\a.\b_1.1+1.\b_2.\a^*+\a.\b_3\a^*$ be as in 5.3(c).
Then $\che\x=1.{}^w\b_0^*.1+\a.{}^w\b_2^*.1+1.{}^w\b_1^*.\a^*+\a.{}^w\b_3^*.\a^*$ and we must show that
$${}^w\b_2^*-{}^w\b_1^*+\a\i{}^w\b_0^*-\a{}^w\b_3^*={}^{sw}(\b_1-\b_2+\a\i\b_0-\a\b_3)^*.$$
Since $\b_i\in N^1$ we have ${}^{s^*}\b_i^*=-b_i^*$; we have also ${}^{s^*}\a^*=-\a^*$ and $sw=ws^*$. Thus 
$$\align&{}^{sw}(\b_1-\b_2+\a\i\b_0-\a\b_3)^*
={}^w({}^{s^*}\b_1^*-{}^{s^*}\b_2^*-(\a^*)\i{}^{s^*}\b_0^*+\a^*{}^{s^*}\b_3^*)\\&
={}^w(-\b_1^*+\b_2^*+(\a^*)\i\b_0^*-\a^*\b_3^*)\endalign$$
as desired.

We prove (d). In this case we must show that
$$\a{}^w\b_2^*-\a{}^w\b_1^*-{}^w\b_0^*+\a^2{}^w\b_3^*={}^{sw}(\a\b_1-\a\b_2-\b_0+\a^2\b_3)^*$$
for $\b_i\in N_0$. We have ${}^{s^*}\b_i^*=b_i^*$; we have also ${}^{s^*}\a^*=-\a^*$ and $sw=ws^*$. Thus 
$$\align&{}^{sw}(\a\b_1-\a\b_2-\b_0+\a^2\b_3)^*\\&
={}^w(-\a^*{}^{s^*}\b_1^*+\a^*{}^{s^*}\b_2^*-{}^{s^*}\b_0^*+\a^{*2}{}^{s^*}\b_3^*)\&
={}^w(-\a^*\b_1^*+\a^*\b_2^*-\b_0^*-\a^{*2}\b_3^*)\endalign$$
as desired.

\subhead 5.14\endsubhead
The involution $f\m f^!$, $N^M@>>>N^M$, given by $f^!(m)=\ff_w(f(\ph(m)))$ induces an involution
$\Th:\un{N^M}@>>>\un{N^M}$ (see 5.1) and also an involution $\un{\un{N^M}}@>>>\un{\un{N^M}}$ denoted again 
by $\Th$. The involution $f\m f^!$, $N'{}^M@>>>N'{}^M$, given by $f^!(m)=\ff_{sw}(f(\ph(m))$ induces an 
involution $\Th':\un{N'{}^M}@>>>\un{N'{}^M}$ (see 5.1) and also an involution
$\un{\un{N'{}^M}}@>>>\un{\un{N'{}^M}}$ denoted again by $\Th'$. Using that ${}^{w\i}\a=\a^*$ and 
${}^{(sw)\i}\a=-\a^*$, we see that:

(a) {\it $\Th:\un{\un{N^M}}@>>>\un{\un{N^M}}$ is $\RR[\a^*]/(\a^{*2})$-linear;
$\Th':\un{\un{N'{}^M}}@>>>\un{\un{N'{}^M}}$ is only $\RR$-linear and satisfies
$\Th'(f\a^*)=-\Th'(f)\a^*$ for $f\in\un{\un{N'{}^M}}$.}
\nl
Note that $\un{\un{N^M}}$, $\un{\un{N'{}^M}}$ are free right $\RR[\a^*]/(\a^{*2})$-modules. Hence for any 
$i$ we have exact sequences of $\RR$-vector spaces
$$0@>>>\un{N^M}_{i-2}@>c>>\un{\un{N^M}}_i@>d>>\un{N^M}_i@>>>0,$$
$$0@>>>\un{N'{}^M}_{i-2}@>c'>>\un{\un{N'{}^M}}_i@>d'>>\un{N'{}^M}_i@>>>0,$$
where $c,c'$ are induced by right multiplication by $\a^*$ and we have
$d\Th=\Th d$, $d'\Th'=\Th'd$, $c\Th=\Th c$, $c'\Th'=-\Th'c'$. (We use (a).)
It follows that for any $i\in\ZZ$ we have
$$\tr(\Th,\un{\un{N^M}}_i)=\tr(\Th,\un{N^M}_i)+\tr(\Th,\un{N^M}_{i-2}),\tag b$$
$$\tr(\Th',\un{\un{N'{}^M}}_i)=\tr(\Th,\un{N'{}^M}_i)-\tr(\Th,\un{N'{}^M}_{i-2}).\tag c$$
Using the isomorphisms in 5.12 and the identities in 5.13 we see that for any $i\in\ZZ$ we have
$$\tr(\Ph,\un{\un{(X^0)^M}}_i)=\tr(\Th,\un{\un{N^M}}_{i-4}),\qua
\tr(\Ph,\un{\un{(X^1)^M}}_i)=\tr(\Th,\un{\un{N^M}}_{i-2}),$$
$$\tr(\Ph,\un{\un{((R.N.R/X)^0)^M}}_i)=\tr(\Th',\un{\un{N'{}^M}}_{i-2}),$$
$$\tr(\Ph,\un{\un{((R.N.R/X)^1)^M}}_i)=\tr(\Th',\un{\un{N'{}^M}}_i).$$
Introducing this into 5.11(a) we deduce
$$\tr(\Ph,\un{\un{(R.N.R)^M}}_i)=\tr(\Th,\un{\un{N^M}}_{i-4})+\tr(\Th,\un{\un{N^M}}_{i-2})
+\tr(\Th',\un{\un{N'{}^M}}_{i-2})+\tr(\Th',\un{\un{N'{}^M}}_i),$$
from which (taking into account (b),(c) and 5.5.(c)) we deduce
$$\align&\tr(\Ps,\un{P}_{i-4})+\tr(\Ps,\un{P}_{i-6})\\&
=\tr(\Th,\un{N^M}_{i-4})+\tr(\Th,\un{N^M}_{i-6})
+\tr(\Th,\un{N^M}_{i-2})+\tr(\Th,\un{N^M}_{i-4})\\&
+\tr(\Th',\un{N'{}^M}_{i-2})-\tr(\Th',\un{N'{}^M}_{i-4})
+\tr(\Th',\un{N'{}^M}_i)-\tr(\Th',\un{N'{}^M}_{i-2}).\endalign$$
We multiply this equality by $v^{-i}$ and sum over all $i$. We get
$$\align&\sum_i\tr(\Ps,\un{P}_{i-4})v^{-i}+\sum_i\tr(\Ps,\un{P}_{i-6})v^{-i}\\&
=\sum_i\tr(\Th,\un{N^M}_{i-4})v^{-i}+\sum_i\tr(\Th,\un{N^M}_{i-6})v^{-i}
+\sum_i\tr(\Th,\un{N^M}_{i-2})v^{-i}\\&+\sum_i\tr(\Th,\un{N^M}_{i-4})v^{-i}
+\sum_i\tr(\Th',\un{N'{}^M}_{i-2})v^{-i}-\sum_i\tr(\Th',\un{N'{}^M}_{i-4})v^{-i}\\&
+\sum_i\tr(\Th',\un{N'{}^M}_i)v^{-i}-\sum_i\tr(\Th',\un{N'{}^M}_{i-2})v^{-i},\endalign$$
that is (using also 5.3(a)):
$$\align&\e^w(R.M.R,\ph')v^{-4}+\e^w(R.M.R,\ph')v^{-6}\\&
=\e^w(M,\ph)v^{-4}+\e^w(M,\ph)v^{-6}+\e^w(M,\ph)v^{-2}+
\e^w(M,\ph)v^{-4}\\&+\e^{sw}(M,\ph)v^{-2}-\e^{sw}(M,\ph)v^{-4}+\e^{sw}(M,\ph)-\e^{sw}(M,\ph)v^{-2},
\endalign$$
where $\ph'$ is as in 5.1. We divide both sides by $v^{-4}+v^{-6}$; we obtain
$$\e^w(R.M.R,\ph')=\e^w(M,\ph)+\e^w(M,\ph)v^2-\e^{sw}(M,\ph)v^2+\e^{sw}(M,\ph)v^4.$$
This proves 5.2(d) and (equivalently) 5.2(b). Theorem 5.2 is proved.

\head 6. Applications\endhead
\subhead 6.1\endsubhead
Theorem 6.2 below describes the action of $T_s+1\in\HH$ in the $\HH$-module $\cm$ (see 3.5, 3.6) for a 
fixed $s\in S$. We set
$$\II'=\{z\in\II, l(z)<l(sz)\},\qua \II''=\{z\in\II, l(z)>l(sz)\}.$$ 
$$\II_e=\{z\in\II;sz=zs^*\},\qua \II_n=\{z\in\II;sz\ne zs^*\},$$
$$\II'_e=\II'\cap\II_e,\qua\II'_n=\II'\cap\II_n,\qua\II''_e=\II''\cap\II_e,\qua\II''_n=\II''\cap\II_n.$$ 
We denote by $w\m\tw$ the involution of $\II$ given by $w\m sw$ if $w\in\II_e$ and $w\m sws^*$ if 
$w\in\II_n$. 

\proclaim{Theorem 6.2} In the $\HH$-module $\cm$ the following identities hold for any $z\in\II$:
$$(T_s+1)a_z=(u+1)(a_z+a_{\tz}) \text{ if }z\in\II'_e,$$
$$(T_s+1)a_z=(u^2-u)(a_z+a_{\tz}) \text{ if }z\in\II''_e,$$
$$(T_s+1)a_z=a_z+a_{\tz} \text{ if }z\in\II'_n,$$
$$(T_s+1)a_z=u^2(a_z+a_{\tz}) \text{ if }z\in\II''_n.$$
\endproclaim
(Recall that $u=v^2$.) We define a map $\II@>>>\II'$, $z\m\hat z$ by $z\m z$ if $z\in\II'$ and $z\m\tz$ if 
$z\in\II''$. For any $z\in\II$ we set 
$$(T_s+1)a_z=\sum_{y\in\II}c_{y,z}a_y$$
where $c_{y,z}\in\ca'$. The following equality (for any $(M,\ph)\in\tC_\sh$) is a reformulation of Theorem 
5.2:
$$\align&\sum_{y\in\II'}\sum_{z\in\II}\e^z(M,\ph)c_{y,z}a_y=
\sum_{y\in\II'_n}(\e^{\ty}(M,\ph)v^4+\e^y(M,\ph))a_y \\&
+\sum_{y\in\II'_e}(\e^{\ty}(M,\ph)v^4-\e^{\ty}(M,\ph)v^2+\e^y(M,\ph)v^2+\e^y(M,\ph))a_y.\endalign$$
Taking $(M,\ph)=(B_x,\ph_x)$ (see 3.2) we see that for any $x\in\II$  we have
$$\sum_{z\in\II}\e^z(B_x,\ph_x)c_{y,z}=\e^{\ty}(B_x,\ph_x)v^4+\e^y(B_x,\ph_x)\text{ if }y\in\II'_n$$  
$$\sum_{z\in\II}\e^z(B_x,\ph_x)c_{y,z}=\e^{\ty}(B_x,\ph_x)(v^4-v^2)+\e^y(B_x,\ph_x)(v^2+1)\text{ if } 
y\in\II'_e.$$
Since the functions $z\m [x\m\e^z(B_x,\ph_x)]$ from $\II$ to the set of maps $\II@>>>\ca'$ are linearly 
independent (see 3.4, 3.5) we deduce that for $y\in\II'_n$, $z\in\II$ we have 
$$c_{y,z}=v^4\text{ if }y=\tz; c_{y,z}=1\text{ if }y=z; c_{y,z}=0\text{ if }z\n\{y,\ty\};$$
and for $y\in\II'_e$, $z\in\II$ we have 
$$c_{y,z}=v^4-v^2\text{ if }y=\tz; c_{y,z}=v^2+1\text{ if }y=z; c_{y,z}=0\text{ if }z\n\{y,\ty\}.$$
Thus for any $z\in\II$ we have
$$(T_s+1)a_z=r_za_{\hat z}+\sum_{y\in\II''}c_{y,z}a_y\tag a$$
where $r_z=v^4$ if $z\in\II''_n$, $r_z=1$ if $z\in\II'_n$, $r_z=v^4-v^2$ if $z\in\II''_e$, $r_z=v^2+1$ if 
$z\in\II'_e$.

We apply $(T_s+1)$ to both sides of (a) and we use that $(T_s+1)^2=(u^2+1)(T_s+1)$ in $\HH$. We obtain
$$\align&(u^2+1)r_za_{\hat z}+\sum_{y\in\II''}(u^2+1)c_{y,z}a_y\\&
=r_zr_{\hat z}a_{\hat z}+r_z\sum_{y\in\II''}c_{y,\hat z}a_y+\sum_{y\in\II''}r_yc_{y,z}a_{\ty}
+\sum_{y\in\II'',y'\in\II''}c_{y,z}c_{y',y}a_{y'}\endalign$$
for any $z\in\II$. Taking the coefficients of $a_y$ with $y\in\II'$ in the two sides of the last equality 
we obtain
$$(u^2+1)r_z\d_{\hat z,y}=r_zr_{\hat z}\d_{\hat z,y}+r_{\ty}c_{\ty,z}.$$
We see that if $y\in\II'$ then $c_{\ty,z}=0$ unless $y=\hat z$ in which case we have
$c_{\ty,z}=r_{\ty}\i r_z((u^2+1)-r_{\hat z})$. The theorem follows.

\subhead 6.3\endsubhead
By Theorem 6.2, the $\HH$-module $\cm$ is identified with the $\HH$-module denoted in \cite{\INV, 0.3} by
$\un{M}$ in such a way that to $a_y\in\cm$ corresponds to $a_y\in\un{M}$ in \cite{\INV}. The duality 
functor $M\m D(M)$ \cite{\SO, 5.9} can be used to define a $\ZZ$-linear map
$\bar{}:\ck_\sh(\tC)@>>>\ck_\sh(\tC)$ which satisfies $\ov{v^n\x}=v^{-n}\ov{\x}$ for $\x\in\ck_\sh(\tC)$, 
$n\in\ZZ$, satisfies $\ov{[B_x,\ph_x]}=[B_x,\ph_x]$ for any $x\in\II$ and satisfies 
$\ov{u\i(T_s+1)\x}=u\i(T_s+1)\ov{\x}$ for any $s\in S$ and any $\x\in\ck_\sh(\tC)$.
It follows that the operator $\bar{}:\ck_\sh(\tC)@>>>\ck_\sh(\tC)$ corresponds under the bijection
$\c'$ in 3.5(a),(b) to the operator $\bar{}:\cm@>>>\cm$ given by \cite{\INV, 0.2}.
It follows that for $x\in\II$, $\tA_x=v^{-l(x)}\sum_{y\in\II;y\le w}\tP_{y,x}(u)a_y\in\cm$ 
(see 3.5) is fixed by the operator $\bar{}:\cm@>>>\cm$ in \cite{\INV, 0.2} where 
$\tP_{y,x}$ are as in 3.5(c). Using \cite{\INV, 0.4}, it follows that for 
$x\in\II$ we have $\tA_x=A_x$ (notation of \cite{\INV, 0.4}) and that for $y\in\II,y\le x$, $\tP_{y,x}$
coincides with the polynomial $P^\s_{x,y}$ introduced in \cite{\INV, 0.4}. Using now 3.5(d), we see that for 
$y\in\II$, $y\le x$ and $\d\in\{1,-1\}$, the following holds:
$$(P_{y,x}(u)+\d P^\s_{y,x}(u))/2\in\NN[u].\tag a$$
This proves Conjecture 9.12 in \cite{\INV}. (In the case where $W$ is a Weyl group this was already known
from \cite{\LV}.)

\subhead 6.4\endsubhead
For $x,y\in W$ we have $c_xc_y=\sum_{z\in W}h_{x,y,z}(u)c_z$ where $h_{x,y,z}(u)\in\NN[u,u\i]$.
Hence for $z,w\in W$ we have $c_zc_wc_{z{*-1}}=\sum_{w'\in W}\ti h_{z,w,w'}(u)c_w$ where
$$\ti h_{z,w,w'}(u)=\sum_{z'\in W}h_{z,w,z'}(u)h_{z',z^{*-1},w'}(u)=\sum_{z'\in W}h_{w,z^{*-1},z'}(u)
h_{z,z',w'}(u).
$$
For $z\in W,w\in\II$ we write $c_zA_w=\sum_{w'\in\II}b_{z,w,w'}(v)A_{w'}$ where $b_{z,w,w'}(v)\in\ca'$. 
For $z\in W$, $w,w'\in\II$ and $\d\in\{1,-1\}$ the following holds:
$$(\ti h_{z,w,w'}(u)+\d b_{z,w,w'}(u))/2\in\NN[u,u\i].\tag a$$
The proof is analogous to that of 6.3(a).
(In the case where $W$ is a Weyl group and $*=1$ this was stated in \cite{\LV, 5.1}.) In particular,

\head 7. The $\HH$-module $\cm_c$\endhead
\subhead 7.1\endsubhead
Let $\le_L,\le_{LR}$ be the preorders on $W$ defined as in \cite{KL}; let $\si_L,\si_{LR}$ be 
the associated equivalence relations on $W$. In this section we fix an equivalence class $c$ for $\si_{LR}$ 
that is, a two-sided cell of $W$. For $w\in W$ we write $w\le_{LR}c$ if $w\le_{LR}w'$ for some $w'\in c$; we
write $w<_{LR}c$ if $w\le_{LR}c$ and $w\n c$. Let $\cm_{\le c}$ (resp. $\cm_{<c}$) be the $\ca'$-submodule 
of $\cm$ generated by the elements $A_x$ with $x\in\II$ such that $x\le_{LR}c$ (resp. $x<_{LR}c$). We show:

(a) {\it $\cm_{\le c}$ is an $\HH$-submodule of $\cm$.}
\nl
With the notation in 6.4 it is enough to show that, if $z\in W$ and $w,w'\in\II$ satisfy 
$b_{z,w,w'}(v)\ne0$ and $w\le_{LR}c$ then $w'\le_{LR}c$. Using 6.4(a) we have $\ti h_{z,w,w'}(u)\ne0$ hence
$\sum_{z'\in W}h_{z,w,z'}(u)h_{z',z\i,w'}(u)\ne0$. It follows that for some $z'\in W$ we have
$h_{z,w,z'}(u)\ne0$ and $h_{z',z\i,w'}(u)\ne0$ hence $w'\le_{LR}z'\le_{LR}w$ and $w'\le_{LR}w$ so that
$w'\le_{LR}c$, as required. 

A similar proof shows:

(b) {\it $\cm_{<c}$ is an $\HH$-submodule of $\cm$.}
\nl
We now define $\cm_c=\cm_{\le c}/\cm_{<c}$. From (a),(b) we see that $\cm_c$ inherits an $\HH$-module
structure from $\cm_{\le c}$. For $x\in\II\cap c$ we denote the image of $A_x\in\cm_{\le c}$ in $\cm_c$
again by $A_x$. Note that $\{A_x;x\in\II\cap c\}$ is an $\ca'$-basis of $\cm_c$.

\subhead 7.2\endsubhead
In the remainder of this paper we assume that $(W,l)$ satisfies the boundedness property in 
\cite{\HEC, 13.2}. (This holds automatically when $W$ is finite or an affine Weyl group, and it probably 
holds in general.) Then the function $\aa:W@>>>\NN$ is defined as in \cite{\HEC, 13.6}.
We recall the following properties:

(i) if $z,z'$ in $W$ satisfy $z\si_{LR}z'$ then $\aa(z)=\aa(z')$;

(ii) if $z,z'$ in $W$ satisfy $z\le_Lz'$ and $\aa(z)=\aa(z')$ then $z\si_Lz'$. 
\nl
(See \cite{\HEC, Ch.14, P4, P9} and \cite{\HEC, Ch.15}; the assumptions 15.1(a),(b) are satisfied by
\cite{\EW}.

In this subsection we fix $s\in S$. For $w\in W$ we set $\e_w=(-1)^{l(w)}$.
Let $y,w\in\II$. As in \cite{\LV, 4.1}, \cite{\INV, 6.1}, we set
$$v^{-l(w)+l(y)}P^\s_{y,w}(v)=\d_{y,w}+\mu'_{y,w}v\i+\mu''_{y,w}v^{-2}\mod v^{-3}\ZZ[v\i]\tag a$$
where $\mu'_{y,w}\in\ZZ,\mu''_{y,w}\in\ZZ$. (When $y\not\le w$ we set $P^\s_{y,w}=0$.) 

As in \cite{\LV, 4.3}, \cite{\INV, 6.2}, 
for any $y,w\in\II$ such that $sy<y<sw>w$, $\e_y=\e_w$, we define $\cm^s_{y,w}\in\ca'$ by
$$\cm^s_{y,w}=\mu''_{y,w}-\sum_{x\in\II;y<x<w,sx<x}\mu'_{y,x}\mu'_{x,w}-\d_{sw,ws^*}\mu'_{y,sw}
+\mu'_{sy,w}\d_{sy,ys^*}.$$

We have the following result.

{\it Let $w\in\II\cap c$. In the $\HH$-module $\cm_c$ we have the following identities:

(b) if $sw<w$, then $c_sA_w=(u+u\i)A_w$,

(c) If $sw>w$, then $c_sA_w=\Xi+\sum_{z\in\II\cap c;sz<z<sw,\e_z=\e_w}\cm^s_{z,w}A_z$;
\nl
where $\Xi$ is given by 

$\Xi=A_{sws^*}$ if $sw\ne ws^*>w$ and $sws^*\in c$,

$\Xi=0$, otherwise.}
\nl
To prove (b),(c) we make use of the formula for $c_sA_w$ given in \cite{\LV, 4.4} (for Weyl groups) and
\cite{\INV, 6.3} in the general case and show that all terms of that formula which involve $(v+v\i)$ belong 
to $\cm_{<c}$ and can therefore be neglected. It is enough to prove the following statements:

(d) If $sw=ws^*>w$ then $sw<_{LR}c$.

(e) If $sw>w$ and $z\in\II$, $\e_z=-\e_w$, $sz<z<sw$, $\mu'_{z,w}\ne0$, then $z<_{LR}c$.
\nl
We prove (d). Since $sw>w$ we have $sw\le_Lw$. If $sw\si_Lw$ then by \cite{\KL, 2.4}, for any $t\in S$ such
that $(sw)t<sw$ we have $wt<w$; in particular, since $sws^*=w<sw$ we have $ws^*<w$, a contradiction. Thus, we
have $sw\not\si_Lw$. From $sw\le_Lw$, $sw\not\si_Lw$, we 
deduce that $sw\not\si_{LR}w$. (If $sw\si_{LR}w$ then $\aa(sw)=\aa(w)$, see (i); from
$sw\le_Lw,\aa(sw)=\aa(w)$ we deduce $sw\si_Lw$ by (ii), a contradiction.) Now (d) follows.

We prove (e). Since $\mu'_{z,w}\ne0$, the coefficient of $v^{l(w)-l(z)-1}$ in $P^\s_{z,w}(v)$ is $\ne0$.
Using 6.3(a) we deduce that the coefficient of $v^{l(w)-l(z)-1}$ in $P_{z,w}(v)$ is $\ne0$. Since 
$sz<z<sw>w$, the last coefficient is known to be equal to $h_{s,w,z}$ (an integer), see \cite{\KL}. Thus we 
have $h_{s,w,z}\ne0$ so that $z\le_Lw$. If $z\si_Lw$ then by \cite{\KL, 2.4}, for any $t\in S$ such that 
$zt<z$ we have $wt<w$; but from $sz<z$, $z\in\II$, we deduce $zs^*<z$ hence $ws^*<w$. From $ws^*<w$ and
$w\in\II$ we deduce $sw<w$, a contradiction. Thus we have $z\not\si_Lw$. From $z\le_Lw$, $z\not\si_Lw$, we 
deduce that $z\not\si_{LR}w$. (If $z\si_{LR}w$ then $\aa(z)=\aa(w)$ by (i); from
$z\le_Lw,\aa(z)=\aa(w)$ we deduce $z\si_Lw$ by (ii), a contradiction.) Now (e) follows.

This completes the proof of (b) and (c).

\subhead 7.3\endsubhead
For $\d\in\{1,-1\}$ let $\cm_c^\d$ be the $\ca'$-submodule of $\cm_c$ generated by 
$\{A_x;x\in\II\cap c,\e_x=\d\}$. From 7.2(b),(c) we see that $\cm_c^\d$ is an $\HH$-submodule of $\cm$.
Clearly, we have $\cm_c=\cm_c^1\op\cm_c^{-1}$ as $\HH$-modules.

\subhead 7.4\endsubhead
The formulas 7.2(b),(c) for the action of $c_s$ in the basis $\{A_x;x\in\II\cap c\}$ of $\cm_c$ are similar 
to those in a $W$-graph (see \cite{\KL}) since the coefficients in the right hand side of 7.2(c) are integer
constants (but unlike the case of $W$-graph these integer constants can in principle depend on $s$). Note 
that the action of left multiplication by $c_s$ in the basis $\{A_x;x\in\II\}$ of $\cm$ is not given by a 
$W$-graph, due to the appearence of terms involving $v+v\i$.

\enddocument